\documentclass[10pt]{amsart}

\usepackage{amssymb}
\usepackage{amsmath}
\newcommand{\ntop}[2]{\genfrac{}{}{0pt}{1}{#1}{#2}}

\let\newpf\proof \let\proof\relax 
\newenvironment{pf}{\newpf[\proofname]}{\qed\endtrivlist}

\def\HS{{\mathrm{HS}}}

\def\DC{{\mathrm{DC}}}

\def\bm{\begin{matrix}}
\def\em{\end{matrix}}

\newcommand{\bt}{\begin{thm}}
\newcommand{\et}{\end{thm}}

\newcommand{\bl}{\begin{lemma}}
\newcommand{\el}{\end{lemma}}

\newcommand{\beq}{\begin{eqnarray}}
\newcommand{\eeq}{\end{eqnarray}}

\def\BBB{{\mathcal{B}}}

\def\RRR{{\mathcal{R}}}

\def\be{\begin{equation}}
\def\ee{\end{equation}}

\def\ba{{\begin{align}}}
\def\ea{{\end{align}}}

\def\H{{\mathbb H}}

\def\0{{\mathbf 0}}

\def\SL{{\mathrm {SL}}}

\def\PSL{{\mathrm {PSL}}}

\newtheorem*{main thm}{Main Theorem}

\newtheorem{thm}{Theorem}[section]
\newtheorem{cor}[thm]{Corollary}

\newtheorem{lemma}[thm]{Lemma}

\theoremstyle{remark}
\newtheorem{rem}{Remark}[section]

\numberwithin{equation}{section}

\def \bn {\hfill \\ \smallskip\noindent}

\theoremstyle{definition}

\newtheorem{definition}{Definition}[section]
\def\proof{\bn {\bf Proof.} }

\def\note#1
{\marginpar
{\tiny $\leftarrow$
\par
\hfuzz=20pt \hbadness=9000 \hyphenpenalty=-100 \exhyphenpenalty=-100
\pretolerance=-1 \tolerance=9999 \doublehyphendemerits=-100000
\finalhyphendemerits=-100000 \baselineskip=6pt
#1}\hfuzz=1pt}

\def\tr{{\text{tr}}}

\newcommand{\dist}{\operatorname{dist}}

\renewcommand{\mod}{\operatorname{mod}}

\newcommand{\C}{{\mathbb C}}
\newcommand{\D}{{\mathbb D}}

\newcommand{\Q}{{\mathbb Q}}
\newcommand{\R}{{\mathbb R}}

\newcommand{\Z}{{\mathbb Z}}

\def\B0{{\bold{0}}}

\catcode`\@=12

\def\Empty{}
\newcommand\oplabel[1]{
  \def\OpArg{#1} \ifx \OpArg\Empty {} \else
  	\label{#1}
  \fi}
		
\newcommand{\comm}[1]{}
\newcommand{\comment}[1]{}

\begin{document}

\title[AC spectrum of the almost Mathieu operator]
{The absolutely continuous spectrum of the almost Mathieu operator}

\author{Artur Avila}

\date{\today}

\address{
CNRS UMR 7599, Laboratoire de Probabilit\'es et Mod\`eles al\'eatoires\\
Universit\'e Pierre et Marie Curie--Bo\^\i te courrier 188\\
75252--Paris Cedex 05, France
}
\curraddr{IMPA, Estrada Dona Castorina 110, Rio de Janeiro, 22460-320,
Brazil}
\email{artur@math.sunysb.edu}

%\address{
%University of California, Irvine, California
%}
%\email{szhitomi@uci.edu}

\begin{abstract}

We prove that the spectrum of the almost Mathieu operator is absolutely
continuous if and only if the coupling is subcritical.  This settles
Problem 6 of Barry Simon's list of Schr\"odinger operator problems
for the twenty-first century.

\end{abstract}

\setcounter{tocdepth}{1}

\maketitle

%\tableofcontents

\section{Introduction}

This work is concerned with the almost Mathieu operator
$H=H_{\lambda,\alpha,\theta}$
defined on $\ell^2(\Z)$ \be \label {defi} (H u)_n=u_{n+1}+u_{n-1}+2 \lambda
\cos(2\pi[\theta+n \alpha])
u_n \ee where $\lambda \neq 0$ is the coupling, $\alpha \in \R \setminus \Q$
is the
frequency and $\theta \in \R$ is the phase.  This is the most studied
quasiperiodic
Schr\"odinger operator, arising naturally as a physical model (see \cite
{L2} for a recent historical account and for the physics background).

We are interested in the decomposition of the spectral measures in atomic
(corresponding to point spectrum),
singular continuous and absolutely continuous parts.  Our main result is the
following.

\begin{main thm}

The spectral measures of the almost Mathieu operator
are absolutely continuous if and only if $|\lambda|<1$.

\end{main thm}

\subsection{Background}

Singularity of the spectral measures for $|\lambda| \geq 1$ had been
previously established (it follows from \cite {LS}, \cite {L},
\cite {AK}).  Thus the Main Theorem reduces to showing absolute continuity
of the spectral measures for $|\lambda|<1$, which is Problem 6 of
Barry Simon's list \cite {S}.

We recall the history of this problem (following \cite {J}).
Aubry-Andr\'e conjectured the following dependence on $\lambda$ of the
nature of the spectral measures:
\begin{enumerate}
\item (Supercritical regime)
For $|\lambda|>1$, spectral measures are pure point,
%\item (Critical regime) For $|\lambda|=1$, spectral measures are singular
%continuous,
\item (Subcritical regime) For $|\lambda|<1$, spectral measures are
absolutely continuous.
\end{enumerate}
A measure-theoretical version of this conjecture was proved by Jitomirskaya
\cite {J}: it holds for almost every $\alpha$ and $\theta$.  Problem 6 of
\cite {S}, which was formulated after (and was likely partially motivated
by) \cite {J} is entirely
about showing that, in the subcritical regime, this
conjecture is in fact true for {\it all} $\alpha$ and $\theta$.

The description of the supercritical regime
turns out to be wrong as stated.  More precisely, for generic
$\alpha$ there can never be point spectrum \cite {G}, \cite {AS},
whatever $\lambda$ and $\theta$ are chosen, and for every $\alpha$ there is
a generic set of $\theta$ for which there is similarly no point
spectrum \cite {JS}.  Thus the result of \cite {J} is essentially the best
possible in the supercritical regime (one can still look for more
optimal conditions on the parameters, which can be sometimes useful for
other purposes, see \cite {AJ1} and \cite {AJ2}).

%It may be that (2) holds as stated.  It was shown to
%hold for almost every $\alpha$ and $\theta$ in \cite {GJLS} (later this was
%improved to for every $\alpha$ and almost every $\theta$ in \cite {AK}), but
%remains an open problem in general.

There was some hope that the description of
the subcritical regime was actually correct as stated, since the work
of Last \cite {L1}, Gesztesy-Simon \cite {GS} (see also
Last-Simon \cite {LS})
%and Jitomirskaya \cite {J}
established that there are absolutely continuous components (of some spectral
measures) for every $\alpha$ and $\theta$ (belief in the conjecture was
however not unanimous, due to lack of any further evidence for generic
$\alpha$).

Two key advances happened recently.
In \cite {AJ2}, the problem was settled for almost every
$\alpha$ and every $\theta$,
and soon later, in \cite {AD} it was settled for every
$\alpha$ (to be precise, for every $\alpha$ that can not be dealt with by
\cite {J}, \cite {AJ1}) and almost every $\theta$.
Those two results are based on quite independent methods, a ramification of
the fact that, as usually happens in quasiperiodic problems,
several aspects of the operators do depend qualitatively on the arithmetics.

\subsection{Outline}

%In order to establish the complete conjecture, we will need to  on
%the foundation provided by both \cite {AJ2} and \cite {AD}.
%Accordingly, the proof splits into two parts that do not interact.
%The arithmetic properties of $\alpha$, more precisely whether it is
%``well approximated by rational numbers'' or not, will decide which of the
%two methods will be applied.

Our proof of the complete conjecture splits into two parts that do not
interact.  The arithmetic properties of $\alpha$, more precisely whether
it is ``well approximated by rational numbers'' or not, will decide which
of the two methods will be applied.

Let $p_n/q_n$ be the continued fraction approximants to
$\alpha$ and let
\be \label {alphaqn}
\beta=\beta(\alpha)=\limsup_{n \to \infty}
\frac {\ln q_{n+1}} {q_n}.
\ee
For our problem, the key distinction is whether $\beta=0$ (the
{\it subexponential regime}) or $\beta>0$ (the {\it exponential regime}).

\subsubsection{The subexponential regime}

In \cite {E}, Eliasson introduced a sophisticated KAM scheme that allowed
him to study the entire spectrum of one-dimensional quasiperiodic
Schr\"odinger operators in the {\it perturbative regime}.  Applied to the
almost Mathieu operator, his results imply that, for frequencies
satisfying the usual Diophantine condition $\alpha \in \DC$, that is
$\ln q_{n+1}=O(\ln q_n)$, and for
$|\lambda|$ sufficiently small (depending on $\alpha$), the spectral
measures of the almost Mathieu operator are absolutely continuous.

In \cite {AJ2}, a non-perturbative method was introduced that, when
applied to the almost Mathieu operator, gives sharp estimates through the
whole subcritical regime for $\alpha \in \DC$.  Absolute continuity of the
spectral measures was then concluded by showing that, after an appropriate
``change of coordinates'', the smallness requirements of the
KAM scheme of Eliasson were satisfied.

In order to extend the conclusions of \cite {AJ2} to the subexponential
regime $\ln q_{n+1}=o(q_n)$,
there are two main difficulties.  The first is that some
key estimates of \cite {AJ2} break down in this setting (essentially for not
achieving exponential decay of Fourier coefficients which is needed to
address the entire $\beta=0$ regime).
The second is that the ``easy path'' consisting of
reducing to a KAM scheme is out of reach.  Indeed, the expected
limit of the KAM method is the Brjuno condition $\sum \frac {\ln
q_{n+1}} {q_n}<\infty$ on $\alpha$, but this is still stronger than
$\ln q_{n+1}=o(q_n)$.  Thus a novel, more robust, approach to
absolute continuity of the spectral measures will need to be implemented.

We notice that the discussion in the subexponential regime yields
significant information which goes beyond the absolute continuity of the
spectral measures (see for instance Remark \ref {dimension} for an example),
and can also be applied to the more general context considered in \cite {AJ2}
(see \S \ref {gen} and \S \ref {generalization}).

%Let us point out that our estimates can be used to get much more
%information on the subexponential regime.  Some of it generalizes
%estimates of \cite {AJ2}, but we get also some
%results which are new even in the Diophantine regime.  Such results will be
%discussed in the appendix.

\subsubsection{The exponential regime}

In the exponential regime, our approach will be to show that each
exponentially close rational approximation $p_n/q_n$
gives a lower bound on the mass
of the absolutely continuous component of a spectral measure, and that
this lower bound converges to the total mass of the spectral measure.

In \cite {AD}, this approach was used to prove absolute continuity of the
integrated density of states, which is the average of the
spectral measures over different $\theta$ (absolute continuity of the
spectral measures for almost every $\theta$ is obtained as a consequence of
this result, by applying \cite {BJ} and \cite {K}).
The key point of \cite
{AD} was to compare averages of the spectral measures (restricted to a large
part of the spectrum) over long sequences
$\{\theta+j q_n\alpha\}_{j=0}^{b_n-1}$ with the corresponding objects for
the periodic operator obtained by replacing $\alpha$ with $p_n/q_n$.  
In such approach, we clearly lose control of individual phases, and one can
not hope to recover a result for every phase by an abstract scheme such as
Kotani's.

Here we will describe a key novel mechanism of ``cancellation'' among
different phases (which we hope will find wider applicability).  We show
that an abnormally small (compared with the
total mass) absolutely continuous component for any $\theta$ implies the
existence of an abnormally large absolutely continuous component for some
$\theta+j q_n \alpha$.  The latter possibility giving a contradiction,
we conclude that all spectral measures have approximately the correct size.

\begin{rem}

Let us mention that the description of the critical regime at this point is
quite accurate but not complete.  One conjectures
(it is explicit in \cite {J}) that for $|\lambda|=1$, for every $\alpha$ and
$\theta$ the spectral measures are singular continuous.
This is proved for every $\alpha$ in the exponential regime and every
$\theta$ (Gordon's Lemma, \cite {G}, \cite {AS}),
almost every $\alpha$ and $\theta$ (\cite {GJLS}),
and it is currently known to hold for every $\alpha$ and
almost every $\theta$ (\cite {AK}).
See also \cite {A} for a recent discussion
including further evidence for the conjecture.

\end{rem}

\subsection{Quasiperiodic Schr\"odinger operators} \label {gen}

Though this work is dedicated to the almost Mathieu operator, several of our
techniques apply to a more general class of
{\it quasiperiodic Schr\"odinger operators}, where the cosine in
(\ref {defi}) is replaced by an arbitrary real analytic periodic
function of the circle.
Though for such operators, we no longer have a sharp phase transition, it
still makes sense to focus on the characterization and understanding
of the regions in the parameters space (including the energy) exhibiting
subcritical, critical and supercritical behavior.

Much recent work has been
developed regarding the ``supercritical'' regime
(characterized, dynamically, by
positive Lyapunov exponent), especially (but not always)
with Diophantine frequencies
(\cite {BG}, \cite {GS1}, \cite {BJ}, \cite {GS2}).  On the other hand,
some ongoing research (mostly still unwritten, but including \cite
{AJ2} and joint works with Fayad and Krikorian) are dedicated towards the
theory of the ``subcritical'' regime
(conjectured in \cite {AJ2} to be
characterized, dynamically, by {\it almost reducibility}),
hopefully without arithmetic restrictions.  Some information on
the intermediate ``critical'' regime can be obtained as a consequence of
these developments (used in conjunction with
renormalization techniques such as \cite {AK}) as well.

Our results in this paper extend partially to this more general setting and
are an important part of the ``subcritical program'' (one of the goals of
which is to prove pure absolutely continuous spectrum in full generality).
It is important to note that the division of the
analysis into subexponential and exponential parts is still relevant,
with new techniques becoming available precisely when
$\beta>0$.  This shows that the extension of the results of \cite {AJ2}
from the usual Diophantine condition to the weaker
subexponential condition plays
an equally important role in the more general setting as well.  Fortunately,
the subexponential part of this work extends fully to the more general
setting, see \S \ref {generalization}.

On the other hand, the analysis of the exponential regime here is still
heavily bound by precise estimates which are only available for the
almost Mathieu operator, and a different analytic approach will need to be
developed to compensate for this.  However, we note that the ``cancellation
technique'' we describe in this paper, which is the only tool available for
deducing pure absolutely continuous spectrum {\it for every phase} in the
case of very Liouvillean frequencies, is robust enough to be applied in the
more general setting once such difficulties are overcome.

{\bf Acknowledgments:} I would like to thank Svetlana Jitomirskaya and David
Damanik for our joint work on \cite {AJ2} and \cite {AD}, which form the
basis on which this work is built on, and for several helpful suggestions
regarding the writing of this paper.  This research was partially conducted
during the period the author served as a Clay Research Fellow.

\section{Preliminaries}

\subsection{Cocycles}

Let $\alpha \in \R$, $A \in C^0(\R/\Z,\SL(2,\C))$.  We call
$(\alpha,A)$ a {\it (complex) cocycle}.
The {\it Lyapunov exponent} is given by the formula
\be
L(\alpha,A)=\lim_{n \to \infty} \frac {1} {n} \int \ln \|A_n(x)\| dx,
\ee
where $A_n$ is defined by
\be
A_n(x)=A(x+(n-1)\alpha) \cdots A(x).
\ee
It turns out (since irrational rotations are uniquely ergodic), that
\be
L(\alpha,A)=\lim_{n \to \infty} \sup_{x \in \R/\Z}
\frac {1} {n} \ln \|A_n(x)\|
\ee
if $\alpha \in \R \setminus \Q$.
We say that $(\alpha,A)$ is {\it uniformly hyperbolic} if there exists a
continuous splitting $\C^2=E^s(x) \oplus E^u(x)$, $x \in \R/\Z$ such that
for some $C>0$, $c>0$, and for every $n \geq 0$,
$\|A_n(x) \cdot w\| \leq C e^{-c n} \|w\|$, $w \in E^s(x)$
and $\|A_n(x)^{-1} \cdot w\| \leq C e^{-cn} \|w\|$, $w \in E^u(x+n \alpha)$.
In this
case, of course $L(\alpha,A)>0$.  We say that $(\alpha,A)$ is {\it bounded}
if $\sup_{n \geq 0} \sup_{x \in \R/\Z} \|A_n(x)\|<\infty$.

Given two cocycles $(\alpha,A^{(1)})$ and $(\alpha,A^{(2)})$,
a {\it (complex) conjugacy}
between them is a continuous $B:\R/\Z \to \SL(2,\C)$ such that
$A^{(2)}(x)=B(x+\alpha)A^{(1)}(x) B(x)^{-1}$ holds.
The Lyapunov exponent is clearly invariant under conjugacies.

We assume now that $(\alpha,A)$ is a {\it real} cocycle, that is,
$A \in C^0(\R/\Z,\SL(2,\R))$.  The notion of real conjugacy (between real
cocycles) is the same as before except that we now ask for $B \in
C^0(\R/\Z,\PSL(2,\R))$.
Equivalently, one looks for $B \in C^0(\R/2\Z,\SL(2,\R))$
satisfying $B(x+1)=\pm B(x)$.  Real conjugacies still preserve the Lyapunov
exponent.

We say that $(\alpha,A)$ is reducible if it is (real) conjugate to a
constant cocycle.

The fundamental group of $\SL(2,\R)$ is isomorphic to $\Z$.  Let
\be
R_\theta=\left (\bm \cos 2 \pi \theta&-\sin 2 \pi \theta\\ \sin 2 \pi
\theta&\cos 2 \pi \theta \em \right ).
\ee
Any $A:\R/\Z \to \SL(2,\R)$ is homotopic to $x \mapsto R_{nx}$ for some
$n \in \Z$ called the degree of $A$ and denoted $\deg A=n$.

Assume now that $A:\R/\Z \to \SL(2,\R)$ is homotopic to the identity.  Then
there exists $\psi:\R/\Z \times \R/\Z \to \R$ and $u:\R/\Z \times \R/\Z \to
\R^+$ such that
\be
A(x) \cdot \left (\bm \cos 2 \pi y \\ \sin 2 \pi y \em \right )=u(x,y)
\left (\bm \cos 2 \pi (y+\psi(x,y)) \\ \sin 2 \pi (y+\psi(x,y)) \em \right
).
\ee
The function $\psi$ is called a {\it lift} of $A$.  Let $\mu$ be any
probability on $\R/\Z \times \R/\Z$ which is invariant by the continuous
map $T:(x,y) \mapsto (x+\alpha,y+\psi(x,y))$, projecting over Lebesgue
measure on the first coordinate (for instance, take $\mu$ as any
accumulation point of $\frac {1} {n} \sum_{k=0}^{n-1} T_*^k \nu$ where
$\nu$ is Lebesgue measure on $\R/\Z \times \R/\Z$).  Then the number
\be
\rho(\alpha,A)=\int \psi d\mu\ \mod\, \Z
\ee
does not depend on the choices of $\psi$ and $\mu$, and is called the
{\it fibered rotation number} of
$(\alpha,A)$, see \cite {JM} and \cite {H}.

The fibered rotation number is invariant under real conjugacies which are
homotopic to the identity.  In general, if $(\alpha,A^{(1)})$ and
$(\alpha,A^{(2)})$ are real conjugate,
$B(x+\alpha)A^{(2)}(x)B(x)^{-1}=A^{(1)}(x)$, and $B:\R/2\Z \to \SL(2,\R)$
has degree $k$ (that is, it is homotopic to $x \mapsto R_{kx/2}$) then
$\rho(\alpha,A^{(1)})=\rho(\alpha,A^{(2)})+k\alpha/2$.

\subsection{$\SL(2,\R)$ action}

Recall the usual action of $\SL(2,\C)$ on $\overline \C$, $\left ( \bm
a & b \\ c & d \em \right ) \cdot z=\frac {az+b} {cz+d}$.

In the following we restrict to matrices $A \in \SL(2,\R)$.
Such matrices preserve $\H=\{z \in \C,\, \Im z>0\}$. The Hilbert-Schmidt
norm of $A=\left
(\bm a & b \\ c & d \em \right )$ is $\|A\|_\HS=(a^2+b^2+c^2+d^2)^{1/2}$. 
Let
$\phi(z)=\frac {1+|z|^2} {2 \Im z}$ for $z \in \H$.  Then $\|A\|_\HS^2=2
\phi(A \cdot i)$.

One easily checks that $\|R_\theta
A\|_\HS=\|AR_\theta\|_\HS=\|A\|_\HS$, so $\phi(R_\theta z)=\phi(z)$.

We notice that $\phi(z) \geq 1$, $\phi(i)=1$ and $|\ln \phi(z)-\ln \phi(w)|
\leq \dist_\H(z,w)$ where $\dist_\H$ is the hyperbolic metric on $\H$,
normalized so that $\dist_\H(a i , i)=|\ln a|$ for $a > 0$.

\subsection{Almost Mathieu operator}

We consider now almost Mathieu operators
$\{H_{\lambda,\alpha,\theta}\}_{\theta \in \R}$.  The definition is the same
as in the introduction, though we will allow $\alpha$ to be a rational
number $p/q$.
The spectrum $\Sigma=\Sigma_{\lambda,\alpha,\theta}$ does not depend on
$\theta$ for $\alpha \in \R \setminus \Q$.  We let $\Sigma_{\lambda,\alpha}$
be this $\theta$-independent set for $\alpha \in \R \setminus \Q$, and we
let $\Sigma_{\lambda,p/q}=\cup_\theta \Sigma_{\lambda,p/q,\theta}$ in the
rational case.
It is the set of $E$ such that $(\alpha,S_{\lambda,E})$ is not uniformly
hyperbolic, with $S_{\lambda,E}$ given by
\be
S_{\lambda,E}(x)=
\left (\bm E-2\lambda \cos 2 \pi x&-1\\1&0 \em \right ),
\ee

The Lyapunov exponent is defined by
$L_{\lambda,\alpha}(E)=L(\alpha,S_{\lambda,E})$.

\begin{thm} [\cite {BJ}, Corollary 2] \label {L}

For every $\alpha \in \R \setminus \Q$, $\lambda \in \R$, $E \in
\Sigma_{\lambda,\alpha}$, we have $L_{\lambda,\alpha}(E)=\max\{\ln
|\lambda|,0\}$.

\end{thm}

\subsubsection{Classical Aubry duality} \label {clas}

Let $\hat H_{\lambda,\alpha,\theta}=\lambda H_{\lambda^{-1},\alpha,\theta}$. 
If $\alpha \in \R \setminus \Q$ then (see \cite {GJLS}) the spectrum of $\hat
H_{\lambda,\alpha,\theta}$ is $\Sigma_{\lambda,\alpha}$.  This reflects an
important symmetry in the theory of the almost Mathieu operators, known as
{\it Aubry duality}.

Classical Aubry duality expresses an algebraic relation between the
families of operators
$\{H_{\lambda,\alpha,\theta}\}_{\theta \in \R}$ and $\{\hat
H_{\lambda,\alpha,\theta}\}_{\theta \in \R}$ which corresponds eigenvectors
with Bloch waves.  In our notation, it is just the computational fact that
if $u:\R/\Z \to \C$ is an
$\ell^2$ function whose Fourier series satisfies
$\hat H_{\lambda,\alpha,\theta} \hat u=E \hat u$, then
$U(x)=\left ( \bm e^{2 \pi i \theta} u(x)\\ u(x-\alpha) \em \right )$
satisfies $S_{\lambda,E}(x) \cdot U(x)=e^{2 \pi i \theta} U(x+\alpha)$.

\subsubsection{The spectral measure} \label {sm}

Fixing a phase $\theta$ and $f \in \ell^2(\Z)$, we let
$\mu^f=\mu^f_{\lambda,\alpha,\theta}$ be
the spectral measure of $H=H_{\lambda,\alpha,\theta}$
corresponding to $f$.  It is defined so that
\be
\langle (H-E)^{-1} f,f \rangle=\int_\R \frac {1} {E'-E} d\mu^f(E')
\ee
holds for $E$ in the resolvent set $\C \setminus \Sigma$.

We set $\mu=\mu^{e_{-1}}+\mu^{e_0}$
(where $\{e_n\}_{n \in \Z}$ is the
canonical basis of $\ell^2(\Z)$).  It is well known that $\{e_{-1},e_0\}$ form
a generating basis of $\ell^2(\Z)$ \cite {CL}, that is, there is no proper
%subset of $l^2(\Z)$
closed subspace of $\ell^2(\Z)$
which is invariant by $H$ and contains $\{e_{-1},e_0\}$.
In particular the support of $\mu$ is
$\Sigma$ and if $\mu$ is absolutely
continuous then all $\mu^f$, $f \in \ell^2$,
are absolutely continuous.  From now on, we restrict our consideration to
$\mu$ which we will call just {\it the} spectral measure.  Notice that, with
our definitions, the spectral measure has total mass $2$.

\subsection{The $m$-functions} \label {m}

The spectral measure $\mu=\mu_{\lambda,\alpha,\theta}$
can be studied through its Borel transform $M=M_{\lambda,\alpha,\theta}$,
\be
M(z)=\int \frac {1} {E'-z} d\mu(E').
\ee
It maps the upper-half plane $\H$ into itself.

For $z \in \H$, there are non-zero solutions $u^\pm$ of $H u^\pm=z u^\pm$
which are $\ell^2$ at $\pm \infty$, defined up to normalization.  Let
\be
m^\pm=\mp \frac {u^\pm_0} {u^\pm_{-1}}.
\ee
Then $m^+$ and $m^-$ map $\H$ holomorphically
into itself.  Moreover, as discussed in \cite
{JL2} (see also \cite {DKL}),
\be
M=\frac {m^+ m^--1} {m^++m^-}.
\ee
The connection with the cocycle acting on $\overline \C$
arises since
\be
S_{\lambda,z}(\theta) \cdot \mp m^\pm(\theta)=\mp
m^\pm(\theta+\alpha).
\ee
Since the holomorphic function
$m^\pm$ maps the upper-half plane into itself, the non-tangential
limits $\lim_{\epsilon \to 0} m^\pm(E+i\epsilon)$ exist
for almost every $E \in \R$, and define a measurable function of $\R$ which
we still denote $m^\pm(E)$.  We will need the following easy consequence of
Kotani Theory (see \cite {R}, Theorem 1.4 for a more general result).

\begin{thm}[see also the more general \cite {R}, Theorem 1.4] \label {s}

For every $\theta$, for almost every
$E$ such that $L_{\lambda,\alpha}(E)=0$, we have
$m^+_{\lambda,\alpha,\theta}(E)=-\overline
{m^-_{\lambda,\alpha,\theta}}(E)$.

\end{thm}

\begin{pf}

It is a key result of Kotani Theory \cite {S2} that the conclusion holds for
almost every $\theta$.  The point here is to extend this to every $\theta$. 
Fix some arbitrary $\theta$, and let $\theta_n \to \theta$ be some sequence
such that the conclusion holds for $\theta_n$.  Let $K=\{E \in \R,\,
L_{\lambda,\alpha}(E)=0\}$.

Let $T:\overline \H \to \overline \D$ be the conformal map
taking $(-1,0,1)$ to $(-1,-i,1)$,
$t^\pm_n(z)=T(m^\pm(\theta_n,T^{-1}(z)))$ and
$t^\pm(z)=T(m^\pm(\theta,T^{-1}(z)))$.  Notice that $t^\pm_n \to
t^\pm$ uniformly on compacts of $\D$.  Let $\eta^\pm_n=t^\pm_n dx$ and
$\eta^\pm=t^\pm dx$, where $dx$ is normalized Lebesgue
measure on $\partial \D$.  By the Poisson formula,
$\eta^\pm_n \to \eta^\pm$ weakly.  Since
$|t^\pm_n| \leq 1$ and $|t^\pm| \leq 1$, we conclude that $\eta^\pm_n|K \to
\eta^\pm|K$.  By the hypothesis on $\theta_n$, $t^+_n=\overline {t^-_n}$
almost everywhere in $K$.  Thus $\eta^+_n|K=\overline
{\eta^-_n}|K$ and passing to the limit, $\eta^+|K=\overline {\eta^-}|K$.  We
conclude that $t^+=-\overline {t^-}$ almost everywhere in $K$, which implies
the result.
\end{pf}

\comm{
\begin{thm}

For almost every $\theta$, for almost every
$E$ such that $L_{\lambda,\alpha}(E)=0$, we have
$m^+_{\lambda,\alpha,\theta}(E)=-\overline {m^-_{\lambda,\alpha,\theta}}(E)
\in \H$.

\end{thm}

Notice that if $E$ is such that $m^+=-\overline {m^-}$ we have
\be \label {LS}
\frac {d} {dE} \mu_{\lambda,\alpha,\theta}(E)=
\frac {1} {\pi} \lim_{\epsilon \to 0} \Im
M(E+i\epsilon)=\frac {1} {\pi} \phi(m(\theta+\alpha)).
\footnote {This
implies that $\frac {d} {dE} \mu \geq \frac {1} {\pi}$ at almost every
$E$ such that $L(E)=0$.  This estimate holds in principle for almost every
$\theta$, but it behaves well with respect to limits, and hence it holds for
every $\theta$ (existence of an absolutely continuous component of the
spectrum for every $\theta$ also follows from \cite {LS}).}
\ee

\begin{cor}

For every $\theta$, for almost every
$E$ such that $L_{\lambda,\alpha}(E)=0$ (hence, for almost every $E \in
\Sigma_{\lambda,\alpha}$ when $|\lambda|<1$), we have
$m^+_{\lambda,\alpha,\theta}(E)=-\overline {m^-_{\lambda,\alpha,\theta}}(E)
\in \H$.

\end{cor}
}

\begin{cor} \label {Y}

Let $0<\lambda<1$, $\alpha \in \R \setminus \Q$.
Then for every $\theta \in \R/\Z$
there exists a measurable function
$m_{\lambda,\alpha,\theta}:\Sigma_{\lambda,\alpha} \to \H$ such that
$S_{\lambda,E}(\theta) \cdot
m_{\lambda,\alpha,\theta}(E)=
m_{\lambda,\alpha,\theta+\alpha}(E)$ and
\be
\frac {d} {dE} \mu_{\lambda,\alpha,\theta}(E)=
\frac{1} {\pi} \phi(m_{\lambda,\alpha,\theta}(E)).
\ee

\end{cor}

\begin{pf}

Let us show that $m=m^+|\Sigma$ has all the properties.  Equivariance is
obvious.  We need to show that $m^+ \in \H$ for almost every $E \in \Sigma$,
and that $\frac {d} {dE} \mu=\frac{1} {\pi} \phi(m^+)$.

First notice that $m^+=-\overline {m^-}$ for almost every $E \in \Sigma$, by
Theorems \ref {L} and \ref {s}.

To show that $m^+ \in \H$ for almost every $E \in \Sigma$, it is enough to
show that the set of $E$ such that $m^+=-\overline {m^-} \in \R \cup
\{\infty\}$ has zero Lebesgue measure.  Otherwise there would be
a positive Lebesgue measure set of $E \in \R$ such that the non-tangential
limit of $m^+$ is $\infty$ or such that the non-tangential limit of
$m^++m^-$ is $0$, both cases giving a contradiction (using that
if the non-tangential limit of either $m^+$ or $m^++m^-$ is
constant in a set of positive Lebesgue measure then $m^+$ or $m^++m^-$ is
constant everywhere).

If $E$ is such that $m^+=-\overline {m^-} \in \H$ we have
\be \label {LS}
\frac {d} {dE} \mu_{\lambda,\alpha,\theta}(E)=
\frac {1} {\pi} \lim_{\epsilon \to 0} \Im
M_{\lambda,\alpha,\theta}(E+i\epsilon)=
\frac {1} {\pi} \phi(m^+_{\lambda,\alpha,\theta}(E)).
\footnote {This
implies that $\frac {d} {dE} \mu \geq \frac {1} {\pi}$.  Notice that even if
the relation $m^+=-\overline {m^-} \in \H$ (almost everywhere with $L=0$)
was only known for a dense subset
of $\theta$, the estimate $\frac {d} {dE} \mu \geq
\frac {1} {\pi}$ (almost everywhere with $L=0$) could then be concluded for
every $\theta$.  Let us point out that we have not used the
$\theta$-independence of
absolutely continuous spectrum obtained in \cite {LS}.
}
\ee
\end{pf}

\comm{
We will consider energies $E+i\epsilon$, $E \in \R$, $\epsilon>0$.
Then there are non-zero solutions $u^\pm$
of $H u^\pm=(E+i\epsilon) u^\pm$ which
are $\ell^2$ at $\pm \infty$, well defined up to normalization.  We define
$m^\pm=m^\pm(x,E+i\epsilon)$ by
\be
m^\pm=\mp \frac {u^\pm_1} {u^\pm_0}.
\ee
Notice that $m^+ \in \H$ and $S_{\lambda,E+i\epsilon}(x+\alpha) \cdot
m^+(x,E+i\epsilon)$.

Let
\be
M(x,E+i\epsilon)=\int \frac {1} {E'-(E+i\epsilon)} d\mu_x(E').
\ee
Notice that $M(E+i\epsilon) \in \H=\{z,\, \Im z>0\}$.
Then, as discussed in \cite {JL2},
\be
M=\frac {m^+ m^--1} {m^++m^-}.
\ee

By Kotani Theory, for almost every $E \in \R$ such that $L(E)=0$, and for
almost every $x \in \R$ (hence, since $0<\lambda<1$, for almost every $E$ in
the spectrum), the non-tangential limits
$m^\pm(x,E)=\lim_{\epsilon \to 0} m^\pm(x,E+i\epsilon)$ exist and
satisfy $m^+(x,E)=-\overline {m^-(x,E)} \in \H$.  It follows that
$M(x,E)=\lim_{\epsilon \to 0} M(x,E+i\epsilon)$ exists and we have
\be
M(x,E)=i \phi(m^+(x,E)).
\ee
Then
\be
\frac {d\mu} {dE}(x)=\lim_{\epsilon \to 0} \Im
M(x,E+i\epsilon)=\phi(m^+(x,E)).
\ee
We want thus to define $\tilde m(x,E)=m^+(x-\alpha,E)$.  The only difficulty

Since the Lyapunov exponent
vanishes in the spectrum, it follows from Kotani theory, \cite {Si2}, that
there exists a full Lebesgue measure set of $x \in \R$, $E \in \Sigma$
such that $\tilde m(x,E) \in \H$ and $A(x) \cdot \tilde m(x,E)=\tilde
m(x+\alpha,E)$, and
\be
\frac {d} {dE} \mu_x(E)=
\frac{1}{2\pi} \int \phi(\tilde m(x,E)) d\theta.
\ee
}

\subsection{Integrated density of states} \label {density}

The integrated density of states is the function
$N=N_{\lambda,\alpha}:\R \to [0,1]$ such that
\be
N(E)=\int_{\R/\Z} \frac {1} {2} \mu_\theta(-\infty,E] d\theta,
\ee
which is a continuous non-decreasing surjective function.
The Thouless formula relates the Lyapunov exponent to
the integrated density of states
\be \label {inte}
L(E)=\int_\R \ln |E'-E| dN(E').
\ee
There is also a relation to the fibered rotation number
\be
N(E)=1-2 \rho(\alpha,S_{\lambda,E})
\ee
where $\rho(\alpha,S_{\lambda,E}) \in [0,1/2]$.

\subsection{Periodic case} \label {periodic}

Let $\alpha=p/q$, and let $A=S_{\lambda,E}$.
The spectrum $\Sigma_{\lambda,p/q,\theta}$ is the set of
all $E$ such that $|\tr A_q(\theta)| \leq 2$, where $A=S_{\lambda,E}$.
The set of $E$ such that $|\tr A_q(\theta)|<2$ is the union of $q$
intervals, and the closure of each interval is called a band.
We order the bands from left to right.  Inside a band, $\tr A_q(\theta)$
is a monotonic function onto $[-2,2]$.

We define $N_{\lambda,p/q,\theta}=\frac {1} {q} \sum_{i=0}^{q-1}
\mu_{\lambda,p/q,\theta+i \alpha}(-\infty,E]$.  Inside the $i$-th band, we
have the formulas
\be \label {n rho}
q N_{\lambda,p/q,\theta}(E)=k-1+(-1)^{q+k-1} 2 \rho(\theta,E)+
\frac {1-(-1)^{q+k-1}} {2}
\ee
where $0<\rho(\theta,E)<1/2$ is such that $\tr A_q(\theta)=2 \cos 2\pi
\rho(\theta,E)$.

In the interior of a band, $\mu_{\lambda,p/q,\theta}$ has a smooth density.
Since $|\tr A_q(\theta)|<2$, there is a well defined fixed point
$m_{\lambda,\alpha,\theta}(E)$ of $A_q(\theta)$ in $\H$.  Then
\be \label {dnp}
\frac {d} {d E} \mu_{\lambda,p/q,\theta}(E)= \frac{1}{\pi}
\phi(m(\theta)).
\ee

\subsection{Bounded eigenfunctions and absolutely continuous spectrum}
\label {b}

%The following is well known (see for instance \cite {JL2}).

\begin{thm} \label {BBB}

Let $\BBB$ be the set of $E \in \R$ such that the cocycle
$(\alpha,S_{\lambda,E})$
is bounded.  Then $\mu_{\lambda,\alpha,\theta}|\BBB$
is absolutely continuous for all $\theta \in \R$.

\end{thm}

This well known result follows from \cite {GP}.  We will actually need an
explicit estimate, contained in
\cite {JL1}, \cite {JL2} (see also \cite {DKL}, page 197 and \cite {S3}).
We give a proof since we found no
reference for the exact statement we need.

\begin{lemma} \label {2000}

We have
$\mu(E-i\epsilon,E+i\epsilon) \leq C \epsilon
\sup_{0 \leq s \leq C \epsilon^{-1}} \|A_s\|_0^2$, where $C>0$ is a
universal constant.

\end{lemma}

\begin{pf}

We have
$\Im M=\frac {\Im m^+ \Im m^-} {|m^++m^-|^2} \left
(\frac {1+|m^+|^2} {\Im
m^+}+\frac {1+|m^-|^2} {\Im m^-} \right )$.  Since
$\Im m^+,\Im m^->0$, $\frac {\Im m^+ \Im m^-} {|m^++m^-|^2}
\leq \frac {1} {2}$ and
\be
\Im M \leq \frac {1} {2} \left (\frac {1+|m^+|^2} {\Im
m^+}+\frac {1+|m^-|^2} {\Im m^-} \right ).
\ee
Clearly
$\Im M(E+i\epsilon) \geq \frac {1} {2 \epsilon} \mu(E-\epsilon,E+\epsilon)$.
so
\be
\frac {1} {2 \epsilon} \mu(E-\epsilon,E+\epsilon) \leq \max \frac
{1+|m^\pm(E+i\epsilon)|^2} {\Im m^\pm(E+i\epsilon)}.
\ee
We want thus to estimate
\be \label {100}
\frac
{1+|m(E+i\epsilon)|^2} {\Im m(E+i\epsilon)} \leq C
\sup_{0 \leq s \leq C \epsilon^{-1}} \|A_s\|_0^2
\ee
for $m=m^+$, $m=m^-$.  By symmetry, we will only consider the case $m=m^+$.
Let $m_\beta=R_{-\beta} \cdot m$.
%where we have made use of the natural action of $\SL(2,\C)$ on
%$\overline \C$,
%$\left (\bm a&b\\c&d \em \right ) \cdot z=\frac {az+b} {cz+d}$.
Those are so-called
$m$-functions for the corresponding half-line problem with appropriate
boundary conditions, see \cite {JL2}, \S 2.
Assume now that $\epsilon^{-1}$ is an integer (the general statement
reduces to this case).
By Proposition 3.9 of \cite {LS} (a consequence of
Theorem 1.1 of \cite {JL1}), such $m$-functions satisfy the bound
\be \label {50}
\Im m_\beta(E+i\epsilon) \leq (5+\sqrt 24)
\sum_{s=0}^{1+\epsilon^{-1}} \|A_s\|_0^2.
\ee

We notice that the quantity $\frac {1+|z|^2} {\Im z}=2\phi(z)$
is invariant under
$R_\beta$.  By choosing $\beta$ appropriately so to maximize $\Im m_\beta$,
$m_\beta$ becomes purely imaginary with $\Im m_\beta \geq 1$, and $\frac
{1+|m|^2} {\Im m} \leq 2 \Im m_\beta$.  Then (\ref {100}) follows
from (\ref {50}).
\end{pf}

\subsection{Corona estimates}

Given a non-zero vector $U \in \C^2$, it is easy to find a matrix with first
column $U$ that belongs to $\SL(2,\C)$.  We just have to solve an equation
of the type $ad-bc=1$, and it is trivial to get estimates on the size of the
solutions.  If $U$ depends holomorphically of a parameter, to obtain a
holomorphic solution of the same problem with good estimates is much more
challenging, and it is related to the famous Corona Theorem of Carleson
\cite {C}.

The Corona Theorem states that if $d \geq 1$ and
$a_i:\D \to \C$, $1 \leq i \leq d$ are bounded holomorphic
functions such that $\max_i |a_i| \geq \epsilon$ pointwise then
there exist bounded holomorphic functions $b_i:\D \to \C$, $1 \leq i \leq d$
such that $\sum a_i b_i=1$.

After the work of Wolff, good estimates on the solutions $b_i$ were
obtained.  For instance, Uchiyama \cite {U} (see Trent \cite {T} for a
published generalization) showed that if
$\delta \leq (\sum |a_i|^2)^{1/2} \leq 1$ pointwise then the
$b_i$ can be chosen such that $(\sum |b_i|^2)^{1/2}
\leq C \delta^{-2} (1-\ln \delta)$, with $C$ independent of $d$.  (Let us
point out that \cite {C} gives an upper bound of the form
$C_d \delta^{-C_d}$ with $C_d$ depending on $d$,
that would be enough for our purposes.)

If instead of functions of the disk one considers functions of an annulus
$\{x \in \C/\Z,\, |\Im x|<a\}$, the conclusion of the Corona Theorem (with
the Uchiyama estimates) is
still valid, and is a consequence of the disk version (because the annulus
is uniformized by the disk and has amenable fundamental group).

The following is an equivalent convenient formulation of the case $d=2$ of
Uchiyama's Theorem for the annulus.

\begin{thm} \label {matrix B}

Let $U:\R/\Z \to \C^2$ be an analytic function.  Assume that $\delta_1 \leq
\|U(x)\| \leq \delta_2^{-1}$ for $|\Im x|<a$.  Then there exists $B:\R/\Z
\to \SL(2,\C)$ with first column $U$ and such that
$\|B\|_a \leq C \delta_1^{-2} \delta_2^{-1} (1-\ln \delta_1 \delta_2)$.

\end{thm}

\section{The subexponential regime}

%We will follow the basic scheme of \cite {AJ2},
%based on quantitative duality.

The approach to the subexponential regime centers around the notion of
quantitative duality, first developed in \cite {AJ2}: almost localization
estimates for the dual operator yield information on the
Fourier series of a ``conjugacy'' to constant.  The almost localization
estimate gives exponential decay away from ``resonant'' sites, but this does
not ensure convergence for all energies (for generic energies, it is
actually divergent).  Still, the estimates yield a good control
of the dynamics.

Though several aspects of the basic scheme of \cite {AJ2} adapt without
difficulties to our weaker conditions, it is clear that
some of the estimates in \cite {AJ2} lose exponential control of the decay
of Fourier coefficients, and hence are too weak to deal with the small
denominators arising in the regime $\beta=0$ (the fight between the decay of
Fourier coefficients and the small denominators happens when we need to
solve the cohomological equation with small error).  This is overcome by the
%but also a tricky
%problem of ``weak repulsion of resonances'', which makes it difficult to get
%good upper bounds on the cocycle.  We will overcome the difficulties
%by the introduction of the notion of ``good resonances'',
%which have slightly better repulsion properties, and the
systematic use of estimates in a definite strip for the truncated
``conjugacies''.  We then need to relate the control of the dynamics
with absolutely continuous spectrum (as described in the introduction, \cite
{AJ2} invokes the KAM approach at this point, which we can not do).
We have good estimates on cocycle
growth in terms of the resonant character of the dual phase, and
bounds on cocycle growth yield upper
bounds on the spectral measures.  We still need estimates connecting the
``parametrization by dual phase'' with the ``parametrization by
energy''\footnote {Each energy usually (almost everywhere)
corresponds to finitely many dual phases, but we have not been able to rule
out (and it is not even heuristically clear that this should be the case,
see footnote 11 of \cite {AJ2})
that for some exceptional set of energies there could be
uncountably many ones.  This is closely related to the
coexistence of both point and singular continuous spectrum for the dual
model.  Happily for us, the exceptional set is very small (with Hausdorff
dimension zero).}, which is done
through a third parametrization, by fibered rotation number.

Another interpretation of the proof is that we give some H\"older
control (in certain scales, we do not actually show full H\"older
continuity here) on the spectral measures, while showing that the
support of the singular part has Hausdorff dimension zero (with good
coverings at the right scales to match the other estimate).

%(the distinction roughly corresponds
%to point or singular continuous spectrum for the dual model).
%exceptionally correspond
%with information collected will
%allow us to establish links between three relevant parametrizations of the
%dynamics: by energy, by fibered rotation number, and by dual phase (this
%last one is not even ``countably to one'' related to the others).  The
%control on cocycle growth, first obtained in terms the dual phase, is
%transfered to the energy parametrization, of Such
%estimates will be the base of the novel argument for absolutely continuous
%spectrum implemented here.
%which will merge all the
%information we will have collected on the dynamics of the
%cocycles.  There will be three parametrizations involved at the same
%time, corresponding to the energy, the fibered rotation number, and the dual
%phase, which will need to be related through estimates.

\subsection{Strong localization estimates}

Let $\alpha \in \R$, $\theta \in \R$, $\epsilon_0>0$.
We say that $k$ is an $\epsilon_0$-resonance if $\|2
\theta-k\alpha\|_{\R/\Z} \leq e^{-|k|\epsilon_0}$ and $\|2
\theta-k\alpha\|_{\R/\Z}=\min_{|j| \leq |k|} \|2 \theta-j\alpha\|_{\R/\Z}$.

\begin{rem}

In particular, there exists always at least one resonance, $0$.  If
$\beta=0$, $\|2 \theta-k\alpha\|_{\R/\Z} \leq e^{-|k|\epsilon_0}$
implies $\|2 \theta-k\alpha\|_{\R/\Z}=
\min_{|j| \leq |k|} \|2 \theta-j\alpha\|_{\R/\Z}$ for $k$ large.

\end{rem}

We order the
$\epsilon_0$-resonances $|n_1| \leq |n_2| \leq ...$.  We say that $\theta$
is $\epsilon_0$-resonant if the set of resonances is infinite.

%Let $\hat H_{\lambda,\alpha,\theta}=\lambda H_{\lambda^{-1},\alpha,\theta}$. 
%If $\alpha \in \R \setminus \Q$, then the spectrum of $\hat
%H_{\lambda,\alpha,\theta}$ is $\Sigma_{\lambda,\alpha}$, see \cite {GJLS}.

\begin{definition}

We say that $\{\hat H_{\lambda,\alpha,\theta}\}_{\theta \in \R}$ (see \S
\ref {clas}) satisfies a
strong localization estimate if there exists $C_0>0$, $\epsilon_0>0$,
$\epsilon_1>0$ such
that for every eigenfunction $\hat H \hat u=E \hat u$
satisfying $\hat u_0=1$ and $|\hat u_k| \leq 1+|k|$, and for every $C_0
|n_j|<k<C_0^{-1} |n_{j+1}|$ we have $|\hat u_k| \leq C_0 e^{-\epsilon_1
|k|}$.

\end{definition}

\begin{thm}[\cite {AJ2}, Theorem 5.1] \label {SLE}

If $\beta=0$ and $0<\lambda<1$ then
$\{\hat H_{\lambda,\alpha,\theta}\}_{\theta \in \R}$
satisfies a strong localization estimate.

\end{thm}

\subsection{A generalization} \label {generalization}

This section can be ignored if one is only
interested in the proof of the Main Theorem.

Let $v:\R/\Z \to \R$ be analytic and let
$H=H_{v,\alpha,\theta}:\ell^2(\Z) \to \ell^2(\Z)$ be the
quasiperiodic Schr\"odinger operator given by
$(Hu)_n=u_{n+1}+u_{n-1}+v(\theta+n\alpha) u_n$.  The almost Mathieu
operator corresponds to the special case $v(\theta)=2 \lambda
\cos 2 \pi \theta$ for some $\lambda \neq 0$.

As for the almost Mathieu case, the spectral properties of
$\{H_{v,\alpha,\theta}\}_{\theta \in \R}$ are intimately connected with the
Schr\"odinger cocycles $\{(\alpha,S_{v,E}\}_{E \in \R}$, where
$S_{v,E}(x)=\left (\bm E-v(x) & -1 \\ 1 & 0 \em \right )$, and several key
notions have identical development, including spectral measures \S \ref {sm},
$m$-functions \S \ref {m} (except for Corollary \ref {Y} which needs to be
reformulated), integrated density of states \S \ref {density} and bounded
eigenfunctions \S \ref {b}.

Most importantly, classical Aubry duality (\S \ref {clas})
can be extended to this setting:
the operators $\hat H_{v,\alpha,\theta}$
given by $(\hat H \hat u)_n=\sum \hat v_k
\hat u_{n-k}+2 \cos (2 \pi (\theta+n \alpha)) \hat u_n$, where $v(x)=\sum
\hat v_k e^{2 \pi i k x}$ have the property that if $u:\R/\Z \to \C$
is an $\ell^2$ function such that
$\hat H_{v,\alpha,\theta} \hat u=E \hat u$, then $S_{v,E}(x) \cdot U(x)=e^{2
\pi i \theta} U(x+\alpha)$, where $U(x)=\left (\bm e^{2 \pi i \theta}
u(x)\\ u(x-\alpha) \em \right )$.

Let us say that $v$ is {\it small} if the family
$\{\hat H_{v,\alpha,\theta}\}_{\theta \in \R}$ is almost localized (the
definition of almost localization being the same as in the almost
Mathieu case).  In particular, $v(x)=2 \lambda \cos 2 \pi x$ is small if
$0<\lambda<1$.  In
general, this notation is justified by Theorem 5.1 of \cite {AJ2} which
shows that if $0<\lambda<\lambda_0(v)$ then
$\{\hat H_{\lambda v,\alpha,\theta}\}_{\theta \in \R}$ is
almost localized in the whole subexponential regime.

We will actually prove the following more general result in the
subexponential regime.

\begin{thm}

If $v$ is small and $\beta=0$ then the spectral measures of
$H_{v,\alpha,\theta}$ are absolutely continous.

\end{thm}

All the discussion below applies essentially unchanged to operators
$H_{v,\alpha,\theta}$ with small $v$ and $\beta=0$.  Besides
replacing mentions of $\lambda$ by $v$ and of the bound
$0<\lambda<1$ by the condition that $v$ is
small, all the few places where
modifications are necessary will be explicitly pointed out in a footnote.

\subsection{Localization and reducibility}

Until the end of this section
we fix $0<\lambda<1$,
$\alpha \in \R \setminus \Q$ with $\beta=0$.
For an energy $E \in \Sigma$, it is shown in Theorem 3.3 of \cite {AJ2}
that there exists some $\theta \in \R$ and $\hat u=(\hat u_i)_{i \in \Z}$
such that $\hat H \hat u=E \hat u$, $\hat u_0=1$,
$|\hat u_i| \leq 1$.  Until the end of this section, whenever $E \in
\Sigma$ is fixed, we will choose some arbitrary $\theta$ and $\hat u$ with
those properties, and we will denote $A=S_{\lambda,E}$.

By the strong localization estimate, if $\theta$ is non-resonant then
$\hat u$ is localized, that is, it is
the Fourier series of an analytic function.
\comm{
Let $u(x)=\sum \hat u_k e^{2 \pi i k x}$ and $U(x)=\left
(\bm e^{2 \pi i \theta} u(x)\\u(x-\alpha) \em \right )$.
}
Classical Aubry duality (\S \ref {clas}) yields a connection between
localization and reducibility (see for instance
Theorem 2.5 of \cite {AJ2}\footnote {Their argument only needs the
arithmetical properties of $\alpha$
to solve the cohomological equation $\phi(x+\alpha)-\phi(x)=
b(x)-\int_0^1 b(x) dx$ with $b$ analytic, and this can be always done when
$\beta=0$.}):

\begin{thm} \label {reducible non-resonant}

If $\theta$ is non-resonant then $(\alpha,A)$ is reducible.

\end{thm}

\comm{
Then $A(x) U(x)=e^{2
\pi i \theta} U(x)$.  Let $B(x)$ be the matrix with columns
$U(x)$ and $\overline U(x)$.  Then $\det B$ is a constant, and if it is
non-zero we have $B(x+\alpha)^{-1} A(x) B(x)=\left (\bm e^{2 \pi i
\theta}&0\\0&e^{-2 \pi i \theta} \em \right )$.
If $\det B$ is zero then
$U(x)=z W(x)$ with $W(x)$ a real vector defined up to sign and $|z|=1$. 
Passing possibly to a double cover, $W$ is well defined, and $W$ does not
vanish over $\R/2\Z$.  Take $\Phi(x)$ with columns $W(x)$ and $\frac {1}
{\|W(x)\|^2} R_{1/4} \Phi(x)$.  Then $\Phi(x+\alpha)A(x)\Phi(x)^{-1}=
\left (\bm 1&b(x)\\0&1 \em \right )$.
We can also further conjugate $A$ to a constant parabolic matrix by
solving the cohomological equation $\phi(x+\alpha)-\phi(x)=b(x)-\int_0^1
b(x) dx$ with $\int_0^1 \phi(x) dx=0$ in $\R/\Z$ (notice that $b$ is
defined in $\R/\Z$): then
$\left (\bm 1&-\phi(x+\alpha)\\0&1 \em
\right ) \Phi(x+\alpha) A(x) \Phi(x)^{-1} \left (\bm 1&\phi(x)\\0&1 \em
\right )=\left (\bm 1&\int_0^1 b(x) dx\\0&1 \em \right )$.
}

\subsection{Bounds on growth}

%The discussion in this section is a straightforward variation of \cite
%{AJ2}.  Besides the weaker Diophantine condition, we make estimates in a
%definite strip, but no difficulty appears yet.

%We start with a general estimate.
%
%\begin{thm}[\cite {BJ}, Corollary 2]
%
%We have $L(\alpha,A)=0$.
%
%\end{thm}

The starting information on the cocycle growth
is given by Theorem \ref {L}, that $L(\alpha,A)=0$.\footnote {For the
generalization, one applies Theorem 6.2 of \cite {AJ2}
whose proof is unchanged in the $\beta=0$ regime.}
In our context this means that for any
$\delta>0$ there exists $c_\delta>0$, $C_\delta>0$ such that
\be \label {context}
\sup_{|\Im x|<c_\delta} \|A_k(x)\| \leq C_\delta e^{-\delta k}.
\ee
The constants $c_\delta$ and $C_\delta$ do not depend on $E$, only on
$\lambda$ and $\alpha$.\footnote {In the case of the
almost Mathieu operator it is possible to show
that we can take $c_\delta=-\frac {1} {2 \pi} \ln \lambda$.  For the
generalization, it is possible to show that it is enough to choose
$c_\delta$ such that $v$ holomorphic in a neighborhood of
$\{|\Im x| \leq c_\delta\}$ and
$c_\delta \leq \frac {1} {2\pi} \epsilon_1$ where
$\epsilon_1$ is the one in the strong localization estimate.}
All further constants may
depend on $\alpha$ and $\lambda$ (respectively $v$).
In the following $C$ is big and $c$ is small.

For a bounded analytic function $f$
defined on a strip $\{|\Im z|<\epsilon\}$ we
let $\|f\|_\epsilon=\sup_{|\Im z|<\epsilon} |f(z)|$.  If $f$ is a bounded
continuous function on $\R$, we let $\|f\|_0=\sup_{x \in \R} |f(x)|$.

Our goal in this section is to prove:

\begin{thm} \label {growth polynomial}

We have $\|A_n\|_c \leq C n^C$.

\end{thm}

Given Fourier coefficients $\hat w=(\hat w_k)_{k \in \Z}$ and an interval $I
\subset \Z$, we let $w^I=\sum_{k \in I} \hat w_k e^{2 \pi i k x}$.  The
length of the interval $I=[a,b]$ is $|I|=b-a$.

We will say that a trigonometrical polynomial $p:\R/\Z \to \C$
has essential degree at most
$k$ if its Fourier coefficients outside
an interval $I$ of length $k$ are vanishing.

Let $p_n/q_n$ be the approximants of $\alpha$.  We recall the
basic properties:
\be \label {b1}
\|q_n \alpha\|_{\R/\Z}=\inf_{1 \leq k \leq q_{n+1}-1}
\|k\alpha\|_{\R/\Z},
\ee
\be \label {b2}
1 \geq q_{n+1} \|q_n \alpha\|_{\R/\Z} \geq 1/2.
\ee
The condition $\beta=0$ implies
\be \label {b3}
q_{n+1} \leq e^{o(q_n)}.
\ee

\comm{
\begin{lemma}[Lemma 9.7 of \cite {AJ1}] \label {9.7}

Let $x \in \R$ and let $0 \leq l_0 \leq q_n-1$ be such that $|\sin \pi
(x+l_0 \alpha)|$ is minimal.  Then
\be \label {irrat-1}
-C \ln q_n \leq \sum_{\ntop {l=0} {l \neq l_0}}^{q_n-1}
\ln |\sin \pi (x+l \alpha)|+(q_n-1) \ln 2 \leq C \ln q_n.
\ee

\end{lemma}
}

\begin{thm}[\cite {AJ2}, Theorem 6.1] \label {trig pol}

Let $1 \leq r \leq [q_{n+1}/q_n]$.
If $p$ has essential degree $k=r q_n-1$ and $x_0 \in \R/\Z$ then
\be \label {d1}
\|p\|_0 \leq C q_{n+1}^{C r} \sup_{0 \leq j \leq k}
|p(x_0+j\alpha)|.
\ee

\end{thm}

In particular, under the condition $\beta=0$
\be \label {d2}
\|p\|_0 \leq C e^{o(k)} \sup_{0 \leq j \leq k} |p(x+j\alpha)|.
\ee

\begin{lemma} \label {following resonance}

We have $o(|n_{j+1}|) \geq \ln \|2 \theta-n_j \alpha\|_{\R/\Z} \geq
c |n_j|$.

\end{lemma}

\begin{pf}

This follows immediately from $\beta=0$.
\end{pf}

Choose $C |n_j|<n<C^{-1} |n_{j+1}|$ of the form $n=r q_k-1<q_{k+1}$, let
$I=[-[n/2],n-[n/2]]$
and define $u(x)=u^I(x)$.  Let
$U(x)=\left (\bm e^{2\pi i \theta} u(x)\\u(x-\alpha) \em
\right )$.  Then
\be
A(x) \cdot U(x)-e^{2 \pi i \theta} U(x+\alpha)=e^{4 \pi i \theta}
\left (\bm h(x)\\0 \em \right ),
\ee
where
\be
\hat h_k=\chi_I(k) 2 \cos 2 \pi (\theta+k \alpha) \hat u_k+\sum_{j \in
\{-1,1\}} \chi_I(k-j) \hat u_{k-j},\footnote {For the generalization one has
$\hat h_k= \chi_I(k) 2 \cos 2 \pi (\theta+k \alpha) \hat u_k+\sum
\chi_I(k-j) \hat v_j \hat u_{k-j}$.}
\ee
where $\chi_I$ is the characteristic function of $I$.
Since $\hat H \hat u=E \hat u$, we also have
\be
-\hat h_k=\chi_{\Z \setminus I}(k) 2 \cos 2 \pi (\theta+k \alpha)
\hat u_k+\sum_{j \in \{-1,1\}} \chi_{\Z \setminus I}(k-j)
\hat u_{k-j}.\footnote {For the generalization one has
$-\hat h_k= \chi_{\Z \setminus I}(k) 2 \cos 2 \pi (\theta+k \alpha)
\hat u_k+\sum \chi_{\Z \setminus I}(k-j) \hat v_j \hat u_{k-j}$.}
\ee
The estimates $|\hat u_k|<C e^{-c |k|}$ for $C^{-1} n<|k|<C n$, $|\hat u_k|
\leq 1$ for all $k$ then
imply that $|\hat h_k| \leq C e^{-c n} e^{-c k}$, that is $\|h\|_c
\leq C e^{-c n}$.

In the following, $\delta$ and $\delta_0$
will be suitably small constants (much smaller
than the $c$ that appeared so far).

\begin{thm} \label {lower bound ux}

We have $\inf_{|\Im x|<\delta_0} \|U(x)\| \geq c e^{-\delta n}$.

\end{thm}

\begin{pf}

Otherwise, by (\ref {context}), $|u(x+j\alpha)| \leq c e^{-\delta n/2}$
for some $x$ with $\Im x=t$, $|t|<\delta_0$ and $0 \leq j \leq n$.  Then
$\|u_t\|_0 \leq c e^{-\delta n/5}$ by
Theorem \ref {trig pol}, where $u_t(x)=u(x+t i)$.
This contradicts $\int u_t(x) dx=1$.
\end{pf}

Let $B(x) \in \SL(2,\C)$ be the matrix whose first column $U(x)$ given by
Theorem \ref {matrix B}.
Then
\be
B(x+\alpha)^{-1} A(x) B(x)=\left (\bm e^{2 \pi i \theta}&0\\0&e^{-2 \pi i
\theta} \em \right )+\left (\bm \beta_1(x)&b(x)\\\beta_3(x)&\beta_4(x)
\em \right ).
\ee
where $\|b(x)\|_{\delta_0} \leq C e^{3 \delta n}$, and
$\|\beta_1(x)\|_{\delta_0},\|\beta_3(x)\|_{\delta_0},
\|\beta_4(x)\|_{\delta_0}
\leq C e^{-c n}$.  Taking $\Phi(x)$ the product of $B(x)^{-1}$
and a constant diagonal matrix, $\Phi(x)=D B(x)^{-1}$, where
$D=\left (\bm d&0\\0&d^{-1} \em \right )$, with $d^2=\max
\{\|\beta_3\|^{1/2}_{\delta_0},e^{-c n}\}$, we get
\be
\Phi(x+\alpha) A(x) \Phi(x)^{-1}=
\left (\bm e^{2 \pi i \theta}&0\\0&e^{-2 \pi i   
\theta} \em \right )+Q(x),
\ee
where $\sup_{|\Im x|<\delta_0} \|Q(x)\| \leq C e^{-c n}$ and
$\sup_{|\Im x|<\delta_0} \|\Phi(x)\| \leq C e^{c n}$.  Thus
\be \label {cecn}
\sup_{0 \leq s \leq c e^{c n}} \|A_s\|_{\delta_0} \leq C e^{c n}.
\ee

\noindent {\it Proof of Theorem \ref {growth polynomial}.}
Let $m \geq C$.  By Lemma \ref {following resonance}
we can choose $C \ln m \leq n \leq C \ln m$
so that $C |n_j|<n<C^{-1} |n_{j+1}|$ and $n=r q_k-1<q_{k+1}$ for some $j$
and $k$.  By (\ref {cecn}), $\|A_m\|_c \leq C m^C$.
\qed

\subsection{Triangularization in a definite strip}

\begin{thm} \label {triangula}

Fix some $n=|n_j|$ and let $N=|n_{j+1}|$ if defined, otherwise let
$N=\infty$.  Then there exists $B:\R/\Z \to \SL(2,\C)$
analytic with $\|B\|_c \leq e^{o(n)}$ such that
\be
B(x+\alpha)A(x)B(x)^{-1}=\left (\bm e^{2 \pi i \theta} & 0\\0 & e^{-2
\pi i \theta} \em \right )+\left (\bm \beta_1(x)&b(x)\\
\beta_3(x)&\beta_4(x) \em \right ),
\ee
with $\|\beta_1\|_c,\|\beta_3\|_c,\|\beta_4\|_c
\leq e^{-cN}$ and $\|b\|_c \leq e^{-c n}$.
In particular
\be \label {good}
\|A_s\|_c \leq C e^{o(n)}, \quad 0 \leq s \leq e^{c n}.
\ee

\end{thm}

\begin{pf}

Let $u(x)=u^I(x)$ for $I=[-c N,c N]$.  Let $r q_k>C n_j$ be minimal with $r
q_k-1<q_{k+1}$ and let
$J=[-[r q_k/2],r q_k-1-[r q_k/2]]$.  Define $U(x)$ as before, and define
also
$U^J(x)$.  Then our previous estimate Theorem \ref {lower bound ux} can be
improved to $\inf_{|\Im x|<c} \|U^J(x)\| \geq e^{-o(n)}$.
The estimate is better since we can use Theorem \ref {growth polynomial}
instead of the weaker estimate (\ref {context}).
Since $\|U-U^J\|_c \leq e^{-c n}$, we get
\be \label {qkc}
\inf_{|\Im x|<c} \|U(x)\| \geq e^{-o(n)}.
\ee
Moreover, we have $A(x) \cdot U(x)=e^{2 \pi i \theta}
U(x+\alpha)+\left (\bm h(x)\\0 \em \right )$ with
$\|h\|_c \leq e^{-c N}$.
Taking $\tilde B$ given by Theorem \ref {matrix B}, we get
\be
\tilde B(x+\alpha)A(x)\tilde B(x)^{-1}=
\left (\bm e^{2 \pi i \theta} & 0\\0 & e^{-2
\pi i \theta} \em \right )+\left (\bm \tilde \beta_1(x)&\tilde b(x)\\
\tilde \beta_3(x)&\tilde \beta_4(x) \em \right ),
\ee
with $\|\tilde \beta_1\|_c,\|\tilde \beta_3\|_c,\|\tilde \beta_4\|_c
\leq e^{-cN}$ and $\|\tilde b\|_c \leq e^{o(n)}$.  If $n \leq C$ we are
done, otherwise let $b^{(1)}(x)$ be obtained by truncating the
Fourier series of
$\tilde b$, so that it has the Fourier coefficients with $|k| \leq n-1$.
We solve exactly
\be
W(x+\alpha) \left (\bm e^{2 \pi i \theta}&b^{(1)}(x)\\0&e^{-2 \pi i
\theta} \em \right )W(x)^{-1}=
\left (\bm e^{2 \pi i \theta}&0\\0&e^{-2 \pi i
\theta} \em \right )
\ee
with $W(x)=\left (\bm 1 & -w(x) \\ 0 & 1 \em \right )$, that is
$b^{(1)}(x)-e^{-2 \pi i \theta} w(x+\alpha)+e^{2 \pi i
\theta} w(x)=0$, or in terms of Fourier coefficients,
$\hat w_k=-\hat b_k \frac {e^{-2 \pi i \theta}}
{1-e^{-2 \pi i(2 \theta-k\alpha)}}$.
So we get $\|W\|_c \leq e^{o(n)}$.  Let $B(x)= 
W(x) \tilde B(x)$.  Noticing that
$\|\tilde b-b^{(1)}\|_c
\leq e^{-cn}$, we obtain the estimates on the coefficients of
$B(x+\alpha)A(x)B(x)^{-1}$.  The second statement follows immediately from
the first.
\end{pf}

\subsection{Lower bounds on the integrated density of states}

\begin{thm} \label {epsilon14}

Let $n=|n_j|$ and let $N=|n_{j+1}|$ if defined, otherwise let $N=\infty$.
Let $C e^{-cN} \leq \epsilon \leq e^{-o(n)}$.  Then there exists
$W:\R/\Z \to \SL(2,\C)$ analytic with
$\|W\|_{c n^{-C}} \leq C \epsilon^{-1/4}$ such that
$Q(x)=W(x+\alpha)A(x)W^{-1}$ satisfies
\be
\|Q\|_0 \leq 1+C \epsilon^{1/2}
\ee

\end{thm}

\begin{pf}

Let $B$ be given by Theorem \ref {triangula}.  Let $D=\left (\bm
d&0\\0&d^{-1} \em \right )$ where
$d=\|B\|_c \epsilon^{1/4}$.  Let $W(x)=D B(x)$.  If
$\epsilon \leq e^{-o(n)}$ we have $\|W\|_0 \leq C \epsilon^{-1/4}$.
Moreover
\be
W(x+\alpha)A(x)W(x)^{-1}=\left (\bm e^{2 \pi i \theta} & 0\\0 & e^{-2
\pi i \theta} \em \right )+\left (\bm q_1(x)&q_2(x)\\
q_3(x)&q_4(x) \em \right ),
\ee
with $\|q_1\|_0,\|q_3\|_0,\|q_4\|_0
\leq C \epsilon^{-1/2} e^{-c N}$ and $\|q_2\|_0 \leq C
\epsilon^{1/2}$.
If $\epsilon \geq C e^{-c N}$ then $\|Q\|_0 \leq 1+C \epsilon^{1/2}$.
\end{pf}

\comm{
\begin{cor}

If $B:\R/\Z \to \SL(2,\C)$ is continuous then
$L(\alpha,B) \leq C \|B-A\|_0^{1/2}$.

\end{cor}

\begin{pf}

It is enough to consider the case when $\epsilon=\|B-A\|_0<c$.  Then
$\epsilon$ is in the range specified by Theorem \ref {epsilon14} for some
$n=|n_j|$.  Let $\tilde B(x)=W(x+\alpha) B(x) W(x)^{-1}$.  Then
$\|\tilde B\|_0 \leq \|Q\|_0+\|W\|_0^2 \|B-A\|_0 \leq 1+C
\epsilon^{1/2}$.
Then $L(\alpha,B)=L(\alpha,\tilde B) \leq \ln \|\tilde B\|_0 \leq C
\epsilon^{1/2}$.
\end{pf}
}

\begin{cor} \label {1/2}

The integrated density of states is $1/2$-H\"older.

\end{cor}

\begin{pf}

%This follows from the Thouless formula $L(E)=\int \ln |E-E'| dN(E')$.
By (\ref {inte}), $L(E+i\epsilon) \geq c
(N(E+\epsilon)-N(E-\epsilon))$ for every $\epsilon>0$.  So it is enough to
show that for $0<\epsilon<c$, $L(E+i \epsilon)<C \epsilon^{1/2}$.  The
condition $0<\epsilon<c$ implies that $\epsilon$ belongs to the range
specified by Theorem \ref {epsilon14} for some $n=|n_j|$.
Let $W$ and $Q$ be given by Theorem \ref {epsilon14}.
Then $L(E+i \epsilon)=L(\alpha,\tilde A) \leq \ln
\|\tilde A\|_0$ where $\tilde
A(x)=Q(x)+W(x+\alpha) \left ( \bm i \epsilon & 0 \\ 0 & 0 \em \right )
W(x)^{-1}$.  Clearly $\ln \|\tilde A\|_0 \leq \|Q\|_0-1+\|W\|_0^2 \epsilon
\leq C \epsilon^{1/2}$, so the result follows.
\end{pf}

%Full $1/2$-H\"older continuity is much harder, and we will not pursue it
%here.  We are only interested in a H\"older upper bound since it allows us
%to get a H\"older lower bound.

\begin{lemma} \label {lower bound on N}

If $E \in \Sigma$ then for $0<\epsilon<1$,
$N(E+\epsilon)-N(E-\epsilon) \geq c
\epsilon^{3/2}$.

\end{lemma}

\begin{pf}

Let $\delta=c \epsilon^{3/2}$.
Since $L(E)=0$, by Thouless formula we have
\be
L(E+i\delta)=\int \frac {1} {2} \ln 1+\frac
{\delta^2} {|E-E'|^2} dN(E').
\ee
We split the integral in $I_1=\int_{|E-E'|>1}$,
$I_2=\int_{\epsilon<|E-E'|<1}$,
$I_3=\int_{\epsilon^4<|E-E'|<\epsilon}$ and
$I_4=\int_{|E-E'|<\epsilon^4}$.  We clearly have $I_1 \leq
c^2 \epsilon^3$.  By Corollary \ref {1/2},
it easily follows that $I_4=\sum_{k
\geq 4} \int_{\epsilon^k>|E-E'|>\epsilon^{k+1}} 1+\frac
{\delta^2} {|E-E'|^2} dN(E') \leq C \sum_{k \geq 4}
\epsilon^{k (1/2)} \ln 1+c^2 \epsilon^{1-2k} \leq C \epsilon^{7/4}$.

Using Corollary \ref {1/2}, we can also estimate, with $m=[-\ln \epsilon]$,
\begin{align}
I_2 &\leq \sum_{k=0}^m
\int_{e^{-k-1}}^{e^{-k}} 1+\frac {\delta^2} {|E-E'|^2} dN(E')\\
\nonumber
&\leq C \sum_{k=0}^m c^2 \epsilon^3 e^{2k+2} e^{-k/2}
\leq C c^2 e^{-3m} e^{(3/2) m} \leq C c^2 \delta.
\end{align}

It follows that $I_3 \geq L(E+i\delta)-c \delta$.  It is well known that
$L(E+i\delta) \geq \delta/10$ for $0<\delta<1$ (see Theorem 1.7 of \cite
{DS}).  Thus $I_3 \geq \delta/20$.
Since $I_2 \leq C (N(E+\epsilon)-N(E-\epsilon)) \ln \epsilon^{-1}$,
the result follows.
\end{pf}

\subsection{Real conjugacies}

Again, fix $n=|n_j|$, $N=|n_{j+1}|$ and let $u(x)=u^I(x)$, $I=[-cN,cN]$.
Let $U(x)$ be defined as before,
and let $\tilde U(x)=e^{\pi i n_j x} U(x)$.
Let $\tilde \theta=\theta-n_j \alpha/2$.
Let $B(x)$ be the matrix with columns
$\tilde U(x)$ and $\overline {\tilde U(x)}$.  Let
$L^{-1}=\|2 \theta-n_j \alpha\|_{\R/\Z}$.
Notice that
\be
A(x) \cdot \tilde U(x)=e^{2 \pi i \tilde \theta} \tilde U(x+\alpha)+\left
(\bm h(x)\\0 \em \right )
\ee
with $\|h(x)\|_c \leq e^{-c N}$.  Notice that in the considerations below,
we must pass to a double cover where the dynamics is like $x \mapsto
x+\alpha/2$, but the condition $\beta=0$ is independent of working with
$\alpha$ or $\alpha/2$.

\begin{thm}

We have
\be
\inf_{x \in \R/\Z} |\det B(x)| \geq c L^{-C}.
\ee

\end{thm}

\begin{pf}

Recall the estimate
\be
\inf_{x \in \R/\Z} \|U(x)\| \geq e^{-o(n)}.
\ee
Minimize over
$\lambda \in \C$, $x \in \R/2\Z$ the quantity $\|\overline {\tilde U(x)}-
\lambda \tilde U(x)\|$.  This gives some
$\lambda_0$, $x_0$.  If the result does not hold then
\be \label {2piij}
\|e^{-2 \pi i j \tilde \theta}
\overline {\tilde U(x_0+j\alpha)}-e^{2 \pi i j \tilde
\theta} \lambda_0 \tilde U
(x_0+j\alpha)\| \leq c L^{-C}, \quad 0 \leq j \leq C L^C.
\ee
This implies that
$\|\overline {\tilde U(x_0+j\alpha)}-\lambda_0 \tilde U(x_0+j\alpha)\|
\leq c L^{-c}$ for $0 \leq j \leq c L^{1-c}$, and as before (first truncating
the Fourier series at scale $C \ln L$)
$\sup_{x \in \R/\Z} \|\overline {\tilde U(x)}-\lambda_0 \tilde U(x)\|
\leq c L^{-c}$.  But taking $j=[L/4]$
in (\ref {2piij}), we get $\|i \overline
{\tilde U(x)}+i \lambda_0 \tilde U(x)\| \leq c L^{-c}$, so that
$\|\tilde U(x)\| \leq c L^{-c}$.  This contradicts
$\|\tilde U(x)\|=\|U(x)\| \geq c e^{-o(n)}$.
\end{pf}

Take now $S=\Re \tilde U$, $T=\Im \tilde U$, and let $\tilde W$
be the matrix with columns $S$ and $\pm T$, so to have $\det \tilde W>0$.
Then
\be
A(x) \cdot \tilde W(x)=
\tilde W(x+\alpha) \cdot R_{\mp \tilde \theta}+
O(e^{-c N}), \quad x \in \R/\Z.
\ee
Let $W(x)=|\det B(x)/2|^{-1/2} \tilde W(x)$ so to have $\det W=1$.  Then
\be
A(x) \cdot W(x)=\frac {|\det B(x+\alpha)|^{1/2}} {|\det B(x)|^{1/2}}
W(x+\alpha) \cdot R_{\mp \tilde \theta}+
O(e^{-c N}), \quad x \in \R/\Z.
\ee
Since $\det W(x)=1$, this gives
\be
A(x) \cdot W(x)=W(x+\alpha) \cdot R_{\mp \tilde \theta}+
O(e^{-c N}), \quad x \in \R/\Z.
\ee

Assume first that $n_j$ is even so that
$W(x+1)=W(x)$ and everything is defined in $\R/\Z$.
Letting $\Phi(x)=W(x+\alpha)^{-1} A(x) W(x)$, we get
$\|\rho(\alpha,\Phi)\pm\tilde \theta\|_{\R/\Z} \leq C e^{-c N}$.
Assume now that $n_j$ is odd, so that
$W(x+1)=-W(x)$.  Letting $\Phi(x)=R_{(x+\alpha)/2}
W(x+\alpha)^{-1} A(x) W(x) R_{-x/2}$, we get $\Phi$ defined in $\R/\Z$ with
\be
\Phi(x)=R_{\frac {\alpha} {2}\mp\tilde \theta}+O(e^{-c N}),
\quad |\Im x|<c L^{-C}.
\ee
Then $\|\rho(\alpha,\Phi)-\frac {\alpha} {2}\pm\tilde \theta\|_{\R/\Z} \leq
C e^{-c N}$.

In either case, for some $k$ with $||k|-n| \leq 1$
\be
|\|2\theta-n_j \alpha\|_{\R/\Z}-\|2 \rho(\alpha,\Phi)-k \alpha\|_{\R/\Z}| \leq
\|2\theta-n_j \alpha\|_{\R/\Z}/10.
\ee

To estimate the topological degree of $W$, it is enough to estimate the
degree of $\frac {M(x)} {\|M(x)\|}$ for $M=S$ or for $M=T$.
Notice that $\|\int e^{-\pi i n_j x} (S(x)+iT(x)) dx\| \geq 1$, and select
$M=S$ or $M=T$ so that $\int \|M(x)\| \geq 1/2$.  We have of course $A(x)
\cdot M(x)=M(x+\alpha)+O(e^{-cn})$, $|\Im x|<c$,
which allows us to estimate
\be
\inf_{x \in \R/\Z} \|M(x)\| \geq c e^{-o(n)}.
\ee
as before.
Truncating the Fourier series of $M$ keeping the $|k|<C n$ the
resulting $\tilde M(x)$ is such that $\|\tilde M(x)-M(x)\| \leq
\|M(x)\|/2$, so we just have to estimate the degree of
$\frac {\tilde M(x)} {\|\tilde
M(x)\|}$.  We do this by counting the number of zeroes of the coordinates of
$\tilde M(x)$, and we get $|\deg W| \leq C n$.  Then
\be \label {rhotheta}
|\|2\theta-n_j \alpha\|_{\R/\Z}-\|2 \rho(\alpha,A)-m \alpha\|_{\R/\Z}| \leq
\|2\theta-n_j \alpha\|_{\R/\Z}/10, \text { for some } |m| \leq C n.
\ee
This implies the following result.

\begin{lemma} \label {resonance easy}

If $\theta$ has a resonance $n_j$ then there exists $|m| \leq C |n_j|$ such
that $\|2 \rho(\alpha,A)-m\alpha\|_{\R/\Z}<2 e^{-\epsilon_0 m}$.

\end{lemma}

\subsection{Proof of the Main Theorem in the case $\beta=0$}

%We now give the argument for absolutely continuous spectrum that avoids the
%reduction to a KAM scheme.  It depends on the estimates on cocycle growth,
%and on the relations between $\theta$, $E$ and $N=1-2 \rho$ that we have
%derived.

Let $\BBB$ be the set of $E \in \Sigma$ such that
$(\alpha,S_{\lambda,E})$ is bounded,
and $\RRR$ be the set of $E \in \Sigma$ such that $(\alpha,S_{\lambda,E})$
is reducible.
By Theorem \ref {BBB},
it is enough to prove that for every $\xi \in \R$,
$\mu=\mu_{\lambda,\alpha,\xi}$
is such that $\mu(\Sigma \setminus \BBB)=0$.

Notice that
$\RRR \setminus \BBB$ has only $E$ such that
$(\alpha,S_{v,E})$ is analytically
reducible to parabolic.  It follows that $\RRR
\setminus \BBB$ is countable: indeed for any such $E$ there exists $k \in
\Z$ such that $2\rho(\alpha,S_{\lambda,E})=k \alpha$ in $\R/\Z$.
If $E \in \RRR$,
any non-zero solution $H_{\lambda,\alpha,\xi} u=E u$ satisfies
$\inf_{n \in \Z} |u_n|^2+|u_{n+1}|^2>0$.  In particular there are no
eigenvalues in $\RRR$,
and $\mu(\RRR \setminus \BBB)=0$.  Thus it is enough to prove that
$\mu(\Sigma \setminus \RRR)=0$.

Let $K_m \subset \Sigma$, $m \geq 0$
be the set of $E$ such that there exists $\theta
\in \R$ and a bounded normalized solution
$\hat H_{\lambda,\alpha,\theta} \hat u=E \hat u$ with a
resonance $2^m \leq
|n_j|<2^{m+1}$.  We are going to show that
$\sum \mu(\overline {K_m})<\infty$.  By Theorem \ref {reducible
non-resonant}, $\Sigma \setminus
\RRR \subset \limsup K_m$.  By the Borel-Cantelli Lemma, $\sum \mu(\overline
{K_m})<\infty$ implies that $\mu(\Sigma \setminus \RRR)=0$.

To every $E \in K_m$, let $J_m(E)$ be an open $\epsilon_m=C e^{-c 2^m}$
neighborhood of $E$.  This is chosen so to have $\sup_{0 \leq
s \leq 10 \epsilon_m^{-1}} \|A_s\|_0 \leq e^{o(2^m)}$ by (\ref
{good}).  By Lemma \ref {2000},
\be
\mu(J_m(E)) \leq C e^{o(2^m)} |J_m(E)|,
\ee
where $| \cdot |$ is used for Lebesgue measure.
Take a finite subcover
$\overline {K_m} \subset \cup_{j=0}^r J_m(E_j)$.
Refining this subcover if necessary,
we may assume that every $x \in \R$ is contained in at most $2$ different
$J_m(E_j)$.

By Lemma \ref {lower bound on N}, $|N(J_m(E))| \geq c |J_m(E)|^2$.  By Lemma
\ref {resonance easy}, if $E \in K_m$ then $\|N(E)-k\alpha\|_{\R/\Z}
\leq C e^{-c 2^m}$ for some $|k|<C 2^m$.  This shows that
$N(K_m)$ can be covered by $C 2^m$ intervals $T_s$ of length
$C e^{-c 2^m}$.  Since $|T_s|<C |N(J_m(E))|$
for any $s$, $E \in K_m$,
there are at most $2C+4$ intervals $J_m(E_j)$ such that $N(J_m(E_j))$
intersects $T_s$.  We conclude that there are at most $C 2^m$ intervals
$J_m(E_j)$.  Then
\be
\mu(K_m) \leq \sum_{j=0}^r \mu(J_m(E_j)) \leq C 2^m e^{o(2^m)} e^{-c 2^m},
\ee
which gives $\sum_m \mu(\overline {K_m}) \leq C$.

\begin{rem} \label {dimension}

In fact this argument shows that the set of energies for which the cocycle
is unbounded has Hausdorff dimension zero.  By Theorem \ref {2000}, this set
contains the set of energies where the spectral measures (and the integrated
density of states) are not Lipschitz.

\end{rem}

\section{The exponential regime}

Our approach to the exponential regime centers around the ``parameter
exclusion'' of energies to achieve control on the dynamics by consideration
of rational approximations.  A key difficulty is that our current
understanding does not lead to a uniform (over the phase space) control, but
merely an integrated control.  While the integrated control is
enough to obtain results
almost sure $\theta$, it is certainly dangerous when trying
to obtain results for all $\theta$.

Let us recall the idea of
\cite {AD} to prove absolute continuity of the
integrated density of states.  There it is shown that for a
rational approximant, a large set of the spectrum can be selected where
we have the pointwise estimate
\be
\frac {d} {dE} N_{\lambda,p/q} \leq (1+o(1)) \frac {d} {dE}
N_{\lambda,\alpha},
\ee
as long as $\alpha$ is exponentially close to $p/q$,
which easily implies that the absolutely continuous part of
$dN_{\lambda,\alpha}$ has mass close to $1$.

It seems a very hard problem to obtain pointwise estimates for the spectral
measures themselves, and we are still unable to show that, for
a large set of energies, the dynamics is under control for every
$\theta$, that is, for the moment we can not prevent
\be
\frac {d} {dE} \mu_{\lambda,p/q,\theta} \geq (1+\epsilon) \frac {d} {dE}
\mu_{\lambda,\alpha,\theta}.
\ee
Our approach will be to show that this can not happen for too many energies,
though which small set of energies must be excluded is not explicit
and could depend on $\theta$.
The idea is to show that such bad situation leads to improved estimates
\be
\frac {d} {dE} \mu_{\lambda,p/q,\theta'} \leq (1-\delta) \frac {d} {dE}
\mu_{\lambda,\alpha,\theta'},
\ee
for some other $\theta'$.  Then, integrating on $E$, we conclude that if
for some phase $\theta$ the total mass of the absolutely continuous part is
less than $2-\epsilon$ (we recall that our spectral measure has total
mass $2$)
then we can find another $\theta'$ for which the total mass of
the absolutely continuous part is greater than $2+\delta$, which is clearly
impossible.

\subsection{Proof of the Main Theorem assuming $\beta>0$}

Throughout this section, we fix $\lambda$, $\alpha$ and $\theta$,
and we assume $\beta>0$.

For every $x=\theta+k \alpha$, $k \in \Z$, let
$\tilde m(x,E)=m_{\lambda,\alpha,x}(E)$ (defined for almost every $E \in
\Sigma_{\lambda,\alpha}$), where $m_{\lambda,\alpha,x}$ is 
given by Corollary \ref {Y}.  It is enough to prove that
\be \label {int_Y}
\frac {1} {2 \pi} \int_{\Sigma_{\lambda,\alpha}}
\phi(\tilde m(\theta,E)) dE=1.
\ee

The hypothesis implies that
\be
\left |\alpha - \frac{p}{q} \right | < e^{-(\beta-o(1)) q}
\ee
for arbitrarily large $q$. Fix some $p/q$ with this property and $q$
large.

For a fixed energy $E$, write $A=A^{(\lambda,E)}$,
$A_n=A^{(\lambda,p/q,E)}_n$ and
$\tilde A_n=A^{(\lambda,\alpha,E)}_n$.

Let $c=\min \{\beta/2,-\ln \lambda/2\}$.

Let $X_{\lambda,p/q,\theta}$ be the set of $E$ such that
$\tr A_q(\theta)=2 \cos 2 \pi
\rho(\theta)$ with $1/q^2<\rho(\theta)<1/2-1/q^2$.

The following plays the role of Lemma 3.1 of \cite {AD}.

\begin{lemma} \label {N3.1}

We have $\mu_{\lambda,p/q,\theta}(\Sigma_{\lambda,p/q,\theta} \setminus
X_{\lambda,p/q,\theta}) \leq 4/q$.

\end{lemma}

\begin{pf}

We have $dN_{\lambda,p/q,\theta}=
\frac {1} {q} \sum_{k=0}^{q-1} \mu_{\lambda,p/q,\theta+k\alpha}$ and
$dN_{\lambda,p/q,\theta}(\Sigma_{\lambda,p/q,\theta} \setminus
X_{\lambda,p/q,\theta})=4/q^2$.
\end{pf}

\begin{lemma}[\cite {AD}, Lemma 3.2]

We have
\be
|\Sigma_{\lambda,p/q} \setminus \Sigma_{\lambda,\alpha}| \leq e^{-(c-o(1))
q}.
\ee

\end{lemma}

In particular,
%\footnote {Strictly speaking, our set $X$ is different from
%the set $X$ of \cite {AD}.  However the only thing being used here is that
%$X \subset \Sigma_{\lambda,p/q}$.}
\be
|X_{\lambda,p/q,\theta} \setminus \Sigma_{\lambda,\alpha}|
\leq e^{-(c-o(1)) q}.
\ee

The following is the appropriate version of Lemma 3.3 of \cite {AD} that we
need.  Though the claim is formally different, since our set
$X_{\lambda,p/q,\theta}$ is somewhat
larger than the set $X_{\lambda,p/q}$
given in \cite {AD}, the proof is exactly the same.

If $E$ belongs to the interior of $\Sigma_{\lambda,p/q,\theta}$, let
$m(\theta,E)$ be the fixed point of
$A_q(\theta)$ in $\H$, as in \S \ref {periodic}.

\begin{lemma} \label {o(1)}

We have
\be
\sup_{E \in X_{\lambda,p/q,\theta}} \sup_{x \in \R}
\ln \phi(m(x,E))=o(q).
\ee

\end{lemma}

The analog of Lemma 3.4 of \cite {AD} holds again with same proof.

\begin{lemma} \label {o(q)}

We have $\mu_{\lambda,p/q,\theta}(X_{\lambda,p/q,\theta} \setminus
\Sigma_{\lambda,\alpha})=o(1)$.

\end{lemma}

It follows from Lemmas \ref {N3.1} and \ref {o(q)} that
\be \label {dE}
\int_{X_{\lambda,p/q,\theta} \cap \Sigma_{\lambda,\alpha}} \phi(m(\theta,E)) dE \geq 2 \pi -o(1).
\ee

%The proof of the analogous of Lemma 3.5 of \cite {AD} fails completely.

We will need to estimate how well the cocycle can be compared with the
rational case.
By Lemma \ref {o(1)}, we have $A(\theta+k\alpha)=B(\theta+(k+1) p/q)
R_{\psi(\theta+k p/q)} B(\theta+k p/q)^{-1}$
with $\ln \|B(\theta)\|=o(q)$, by taking
$B(\theta+k p/q) \cdot i=m(\theta+k p/q,E)$. We have
$\prod_{i=q-1}^0 R_{\psi(\theta+i p/q)}=B(\theta)^{-1} A_q(\theta)
B(\theta)= R_{\pm \rho(\theta)}$.
The following is the appropriate version of
Lemma 4.1 of \cite {AD}.  Let $b=[e^{c/10 q}]$.

\begin{lemma} \label {orbit}

For $0 \leq k<b$ we have $\|B(\theta)^{-1} \tilde A_{k q}(\theta)
B(\theta)-R_{\pm k\rho(\theta)}\|=O(e^{-cq/4})$.

\end{lemma}

\begin{pf}

Write
\be
\tilde A_k(\theta)=\prod_{i=k-1}^0
A(\theta+i\alpha)=\prod_{i=k-1}^0
B(\theta+(i+1)p/q) Q_i B(\theta+i p/q)^{-1}.
\ee
Then $\|Q_i-R_{\psi(\theta+ip/q)}\|=O(e^{-(3c/2-o(1))q})$ for $0 \leq i<b q$. 
Thus $
\tilde
A_{kq}(\theta)=B(\theta) Q B(\theta)^{-1}$ where $Q=\prod_{i=kq-1}^0 Q_i$
satisfies
$\|Q-R_{\pm k \rho(\theta)}\|=O(e^{-(c-o(1))q})$.
\end{pf}

We now diverge from \cite {AD}.  The cancellation mechanism will evolve
along the next four lemmas.  The basis is an
equality for periodic elliptic matrices in $\SL(2,\R)$:

\begin{lemma} \label {cancellation}

Let $r/s$ be a rational number which is not an integer multiple of $1/2$,
and let $z_0 \in \H$.  If $B_0 \in \SL(2,\R)$
\be
\frac {1} {s} \sum_{k=0}^{s-1} \phi(B_0 \cdot R_{r k/s} \cdot z_0)=\phi(z_0)
\phi(B_0 \cdot i).
\ee

\end{lemma}

\begin{pf}

Let $Z \in \SL(2,\R)$ be a matrix taking $i$ to $z_0$.  We want to estimate
\be
\frac {1} {s} \sum_{k=0}^{s-1} \|B_0 \cdot R_{r k/s} \cdot Z\|^2_\HS=
\frac {1} {2} \|Z\|_\HS^2 \|B_0\|^2_\HS.
\ee
By considering rotations, this is the same as showing
\be
\frac {1} {s} \sum_{k=0}^{s-1} \|\left ( \bm \nu & 0 \\ 0 & \nu^{-1}
\em \right ) \cdot R_{x+k r/s} \cdot \left ( \bm \nu'
& 0 \\ 0 & \nu'^{-1} \em \right )\|^2_\HS=
\frac {1} {2} \|\left ( \bm \nu' & 0 \\ 0 & \nu'^{-1}\em \right
)\|_\HS^2 \|\left ( \bm \nu & 0 \\ 0 & \nu^{-1} \em \right
)\|^2_\HS,
\ee
for any $x$, $\nu$ and $\nu'$.  A direct computation gives
\begin{align}
&\frac {1} {s} \sum_{k=0}^{s-1} \|\left ( \bm \nu & 0 \\ 0 & \nu^{-1}
\em \right ) \cdot R_{x+k r/s} \cdot \left ( \bm \nu'
& 0 \\ 0 & \nu'^{-1} \em \right )\|^2_\HS=\\
\nonumber
&\frac {1} {s} \sum_{k=0}^{s-1}
\nu^2 \nu'^2+\nu^{-2}
\nu'^{-2}+(\nu^2-\nu^{-2})(\nu'^2-\nu'^{-2}) (\cos^2 2
\pi (x+k r/s)-1)\\
\nonumber
&=\frac {1} {2} (\nu^2+\nu^{-2}) (\nu'^2+\nu'^{-2}),
\end{align}
since
\be
\frac {1} {s} \sum_{k=0}^{s-1} \cos^2 2 \pi(x+k r/s)=\frac {1} {s}
\sum_{k=0}^{s-1} \frac {1+\cos 4 \pi (x+k r/s)} {2}=1/2.
\ee
\end{pf}

Next we obtain an estimate on general elliptic $\SL(2,\R)$ matrices:

\begin{lemma}

For every $\epsilon>0$ there exists $\delta>0$ such that
if $B_0 \in \SL(2,\R)$ with $\|B_0\|<e^{\delta q}$,
$e^{-\delta q}<\rho<1/2-e^{-\delta q}$, $z_0 \in \H$ with
$\phi(z_0)<e^{\delta q}$ and $b_0>e^{\epsilon q}$ then as $q$ grows we have
\be
\frac {1} {b_0}
\sum_{k=0}^{b_0-1} \phi(B_0 \cdot R_{\pm
k \rho} \cdot z_0)>(1-o(1)) \phi(z_0) \phi(B_0 \cdot i).
\ee

\end{lemma}

\begin{pf}

Let $r/s$, $s \leq b_0^{1/2}$ maximal, be an approximant of $\rho$.
Let $z_j=R_{\pm j s \rho} \cdot z_0$, $j \geq 0$.
We can break the sum $\sum_{k=0}^{b_0-1} \phi(B_0 \cdot R_{\pm
k \rho} \cdot z_0)$ into $[b_0/s]$ parts of the form
$\sum_{k=0}^{s-1} \phi(B_0 \cdot R_{\pm
k \rho} \cdot z_j)$, $0 \leq j \leq
[b_0/s]-1$, plus a non-negative remainder.
Since $b_0-s [b_0/s]=o(b_0)$, it is enough to show that
\be
\frac {1} {s}
\sum_{k=0}^{s-1} \phi(B_0 \cdot R_{\pm k \rho} \cdot z_j)>
(1-o(1)) \phi(z_j) \phi(B_0 \cdot i)
\ee
for each $j$ (using that $\phi(z_j)=\phi(z_0)$ for every $j$).
Consider the points $w'_k=B_0 \cdot R_{\pm k r/s} \cdot z_j$, $w_k=
B_0 \cdot R_{\pm k \rho} \cdot z_j$.  Then $\dist_\H(w_k,w'_k)=o(1)$, $0
\leq k \leq s-1$.
It follows that it is enough to show that
\be
\frac {1} {s} \sum_{k=0}^{s-1} \phi(w'_k)>(1-o(1)) \phi(z_j) \phi(B_0 \cdot
i).
\ee
This follows from Lemma \ref {cancellation}.
\end{pf}

We now apply the previous estimates to the cocycle, using that it is well
shadowed by rotations by Lemma \ref {orbit}.

\begin{lemma}[Cancellation along orbits, fixed energy] \label {dynamical}

Let $z \in \H$, $E \in X_{\lambda,p/q,\theta}$.
If $1 \leq \kappa \leq 2$ is such that
\be
|\ln \phi(z)-\ln \phi(m(\theta,E))| \geq \ln \kappa
\ee
then
\be
\frac {1} {b} \sum_{k=0}^{b-1} \phi(\tilde A_{kq}(\theta) \cdot z) \geq
(\frac {1+\kappa^2} {2\kappa}-o(1)) \phi(m(\theta,E)).
\ee

\end{lemma}

\begin{pf}

There are two cases to consider.  If $\phi(z)$ is not $e^{o(q)}$, then by
Lemma \ref {orbit}, $\phi(\tilde A_{kq}(\theta) \cdot z)$
is not $e^{o(q)}$ for $0 \leq k \leq b-1$.  On the other hand,
$\phi(m(\theta,E))=e^{o(q)}$ by Lemma \ref {o(1)}, hence
\be 
\frac {1} {b} \sum_{k=0}^{b-1} \phi(\tilde A_{kq}(\theta) \cdot z) \geq 2
\phi(m(\theta,E)). 
\ee

Assume now that $\phi(z)=e^{o(q)}$.  Set $B_0=B(\theta)$,
$z_0=B(\theta)^{-1} \cdot z$, $\rho=\rho(\theta)$, $b_0=b$ in the previous
lemma.  Then
\be
\frac {1} {b}
\sum_{k=0}^{b-1} \phi(B(\theta) \cdot R_{\pm k \rho(\theta)} B(\theta)^{-1}
\cdot z)>(1-o(1)) \phi(z_0) \phi(B(\theta) \cdot i).
\ee
We have $B(\theta) \cdot i=m(\theta,E)$.
Moreover $\dist_\H(z,m(\theta,E)) \geq
\ln \kappa$ (since $\ln \phi$ is $1$-Lipschitz in the hyperbolic metric), so
we also have $\dist_\H(z_0,i) \geq \ln \kappa$.  It follows that $\phi(z_0)
\geq (1+\kappa^2)/2\kappa$, so we have
\be
\frac {1} {b}
\sum_{k=0}^{b-1} \phi(B(\theta) \cdot R_{\pm k \rho(\theta)} B(\theta)^{-1}
\cdot z)>(1-o(1)) \frac {1+\kappa^2} {2\kappa} \phi(m(\theta,E)).
\ee
We now notice that $\dist_\H(B(\theta) \cdot R_{\pm k \rho(\theta)}
B(\theta)^{-1} \cdot z,\tilde A_{kq} \cdot z)=o(1)$ by Lemma \ref {orbit}
and $\phi(z)=e^{o(q)}$, so
$(1-o(1)) \phi(B(\theta) \cdot R_{k\rho(\theta)} B(\theta)^{-1}
\cdot z) \leq \phi(\tilde A_{kq} \cdot z)$, which gives
\be
\frac {1} {b}
\sum_{k=0}^{b-1} \phi(\tilde A_{kq}
\cdot z)>(1-o(1)) \frac {1+\kappa^2} {2\kappa} \phi(m(\theta,E)).
\ee
\end{pf}

This estimate can be applied to the case $z=\tilde m(\theta)$ and integrated
to yield:

\begin{lemma}[Cancellation along orbits, integrated version]

For every $\epsilon>0$ there exists $\delta>0$ such that if
\be
\int_{X_{\lambda,p/q,\theta} \cap \Sigma_{\lambda,\alpha}}
\phi(\tilde m(\theta,E))<
(1-\epsilon-o(1)) 2\pi
\ee
then
\be
\frac {1} {b} \sum_{k=0}^{b-1} \int_{X_{\lambda,p/q,\theta} \cap
\Sigma_{\lambda,\alpha}}
\phi(\tilde m(\theta+k\alpha,E))>(1+\delta-o(1)) 2 \pi.
\ee

\end{lemma}

\begin{pf}

Let $W \subset X_{\lambda,p/q,\theta} \cap \Sigma_{\lambda,\alpha}$
be the set such that
$\phi(\tilde m(\theta,E))<(1-\epsilon/2) \phi(m(\theta,E))$.  Then by (\ref
{dE})
\be
\int_W \phi(m(\theta,E))>\epsilon \pi-o(1).
\ee
Applying the previous lemma with $z=\tilde m(\theta,E)$ we get
\be
\frac {1} {b} \sum_{k=0}^{b-1} \phi(\tilde
m(\theta+k\alpha,E))>(1+\delta-o(1)) \phi(m(\theta,E)),\quad E \in
X_{\lambda,p/q,\theta} \cap \Sigma_{\lambda,\alpha} \setminus W,
\ee
\be
\frac {1} {b} \sum_{k=0}^{b-1} \phi(\tilde
m(\theta+k\alpha,E))>(1+\delta-o(1)) \phi(m(\theta,E)),\quad E \in W,
\ee
Integrating and using (\ref {dE}) we get
\be
\frac {1} {b} \sum_{k=0}^{b-1} \int_{X_{\lambda,p/q,\theta} \cap
\Sigma_{\lambda,\alpha}}
\phi(\tilde m(\theta+k\alpha,E))>(1-o(1)) 2 \pi+(\delta-o(1))
(\epsilon \pi-o(1)).
\ee
\end{pf}

The conclusion of this lemma being obviously impossible for
$q$ large, we must have
\be
\int_{X_{\lambda,p/q,\theta} \cap \Sigma_{\lambda,\alpha}} \phi(\tilde m(\theta,E))>(1-o(1)) 2 \pi,
\ee
which implies (\ref {int_Y}) as $q$ grows.


\begin{thebibliography}{BKNS}

\bibitem[AA]{AA} Aubry, S.; Andr\'e, G.
Analyticity breaking and
Anderson localization in incommensurate lattices.  Group theoretical methods
in physics (Proc. Eighth Internat. Colloq., Kiryat Anavim, 1979),  pp.
133--164, Ann. Israel Phys. Soc., 3, Hilger, Bristol, 1980.

\bibitem[A]{A} Avila, A.
On point spectrum at critical coupling.  In preparation.

\bibitem[AD]{AD} Avila, A.;Damanik, D.
Absolute continuity of the integrated density of states for the almost
Mathieu operator.  Inv. Math. 172 (2008), 439-453.

\bibitem[AJ1]{AJ1} Avila, A.;Jitomirskaya, S.
The Ten Martini Problem.  Preprint (www.arXiv.org).
To appear in Annals of Math.

\bibitem[AJ2]{AJ2} Avila, A.;Jitomirskaya, S.
Almost localization and almost reducibility.  To appear in Journal of the
European Mathematical Society.

\bibitem[AK]{AK} Avila, A.; Krikorian, R.
Reducibility or non-uniform
hyperbolicity for quasiperiodic Schr\"odinger cocycles.  Annals of Math. 164
(2006), 911-940.

%\bibitem[B1]{B1} H\"older regularity of
%   integrated density of states for the almost Mathieu operator in a
%   perturbative regime. Lett. Math. Phys. 51 (2000), no. 2, 83--118.

%\bibitem[B2]{B2} Bourgain, J. Green's function
%estimates for lattice Schr\"odinger operators and applications. Annals
%of Mathematics Studies, 158. Princeton University Press, Princeton,
%NJ, 2005. x+173 pp.

%\bibitem[AK2]{AK2} Avila, A.; Krikorian, R.
%Quasiperiodic $\SL(2,\R)$ cocycles.  In preparation.

%\bibitem[AvMS]{AvMS} Avron, J.; van Mouche, P. H. M.; Simon, B. On the
%measure of the spectrum for the almost Mathieu operator.  Comm. Math. Phys. 
%132  (1990),  no. 1, 103--118.

\bibitem[AS]{AS} Avron J.;Simon, B.
Singular continuous spectrum for a class of almost periodic
Jacobi matrices.  Bull. Amer. Math. Soc. 6 (1982),
81-85.

\bibitem[BG]{BG} Bourgain, J.; Goldstein, M. On nonperturbative localization
with quasi-periodic potential.  Ann. of Math. (2)  152  (2000),
835--879.

\bibitem[BJ]{BJ} Bourgain, J.; Jitomirskaya, S. Continuity of the Lyapunov
exponent for quasiperiodic operators with analytic potential.
%Dedicated to
%David Ruelle and Yasha Sinai on the occasion of their 65th birthdays.
J. Statist. Phys.  108  (2002), 1203--1218.

%\bibitem[BS]{BS} B\'ellissard, J.; Simon, B.
%Cantor spectrum for the almost Mathieu equation.
%J. Funct. Anal. 48 (1982), no. 3, 408--419.

\bibitem[C]{C} Carleson, Lennart Interpolations by bounded analytic
functions and the corona problem.  Ann. of Math. (2)  76  1962 547--559.

\bibitem[CL]{CL} Carmona, Ren\'e; Lacroix, Jean Spectral theory of random
Schrödinger operators. Probability and its Applications.
Birkh\"auser Boston, Inc., Boston, MA, 1990. xxvi+587 pp.

%\bibitem[CEY]{CEY} Choi, Man Duen; Elliott, George A.; Yui, Noriko Gauss
%polynomials and the rotation algebra.  Invent. Math.  99  (1990),  no. 2,
%225--246.

\bibitem[DKL]{DKL} Damanik, David; Killip, Rowan; Lenz, Daniel Uniform
spectral properties of one-dimensional quasicrystals. III.
$\alpha$-continuity.  Comm. Math. Phys.  212  (2000), 191--204.

\bibitem[DS]{DS} Deift, P.; Simon, B.
Almost periodic Schr\"odinger operators. III. The absolutely continuous
spectrum in one dimension.  Comm. Math. Phys.  90  (1983), 389--411.

\bibitem[E]{E} Eliasson, L. H.
Floquet solutions for the $1$-dimensional quasi-periodic
Schr\"odinger equation.
Comm. Math. Phys. 146 (1992), 447--482.

\bibitem[GS]{GS} Gesztesy, Fritz; Simon, Barry The xi function.
Acta Math. 176 (1996), 49-71.

\bibitem[GP]{GP} Gilbert, D.J., On subordinacy and analysis of the spectrum
of one-dimensional Schr\"odinger operators.  J. Math. Anal. Appl. 128
(1987), 30-56.

\bibitem[GoSc1]{GS1} Goldstein, Michael; Schlag, Wilhelm H\"older
continuity of the integrated density of states for quasi-periodic
Schr\"odinger equations
and averages of shifts of subharmonic functions.  Ann. of Math. (2)  154
(2001), 155--203.

\bibitem[GoSc2]{GS2} Goldstein, Michael; Schlag, Wilhelm
On resonances and the formation of gaps in the spectrum of quasi-periodic
Schr\"odinger equations.  Preprint (www.arXiv.org).

\bibitem[G]{G} Gordon, A.
On the point spectrum of the one-dimensional Schr\"odinger operator.
Usp. Math. Nauk. 31 (1976), 257-258.

\bibitem[GJLS]{GJLS} Gordon, A. Y.; Jitomirskaya, S.; Last, Y.; Simon, B.
Duality and singular continuous spectrum in the almost Mathieu equation. 
Acta Math. 178 (1997), 169--183.

%\bibitem[GS]{GS} Goldstein, M.; Schlag, W.
%Fine properties of the integrated density of states and a quantitative
%separation property of the Dirichlet eigenvalues.  Preprint (www.arXiv.org).

\bibitem[H]{H} Herman, Michael-R. Une m\'ethode pour minorer les
exposants de Lyapounov et quelques exemples montrant le caract\`ere local
d'un th\'eor\`eme
d'Arnol'd et de Moser sur le tore de dimension $2$.
Comment. Math. Helv.  58  (1983), 453--502.

%\bibitem[JK]{JK} Jitomirskaya, S. Ya.; Krasovsky, I. V. Continuity of the
%measure of the spectrum for discrete quasiperiodic operators.  Math. Res.
%Lett.  9  (2002),  no. 4, 413--421.

\bibitem[J]{J} Jitomirskaya, Svetlana Ya. Metal-insulator transition for the
almost Mathieu operator.  Ann. of Math. (2)  150  (1999),
1159--1175.

\bibitem[JL1]{JL1} Jitomirskaya, Svetlana; Last, Yoram Power-law subordinacy
and singular spectra. I. Half-line operators.  Acta Math.  183  (1999),
171--189.

\bibitem[JL2]{JL2} Jitomirskaya, Svetlana; Last, Yoram
Power law subordinacy and singular spectra. II. Line operators.  Comm. Math.
Phys.  211  (2000), 643--658.

\bibitem[JS]{JS} Jitomirskaya, Svetlana; Simon, Barry
Operators with singular continuous spectrum: III. Almost periodic
Schr\"odinger operators.  Comm. Math. Phys. 165 (1994), 201-205.

\bibitem[JM]{JM} Johnson, R.; Moser, J. The rotation number for almost
periodic potentials.  Comm. Math. Phys.  84  (1982), 403--438.

\bibitem[L1]{L} Last, Y.
Zero measure spectrum for the almost Mathieu operator.  Comm.
Math. Phys.  164  (1994), 421--432.

%\bibitem[L2]{L1} Last, Y. Almost everything about the almost Mathieu
%operator. I.  XIth International Congress of Mathematical Physics (Paris,
%1994),  366--372, Internat. Press, Cambridge, MA, 1995.

\bibitem[L2]{L1} Last, Y. A relation between a.c. spectrum of ergodic Jacobi
matrices and the spectra of periodic approximants.  Comm. Math. Phys.  151
(1993), 183--192.

\bibitem[L3]{L2} Last, Y. Spectral theory of Sturm-Liouville operators on
infinite intervals: a review of recent developments.  Sturm-Liouville
theory,  99--120, Birkh\"auser, Basel, 2005.

\bibitem[LS]{LS} Last, Yoram; Simon, Barry
Eigenfunctions, transfer matrices, and absolutely continuous spectrum of
one-dimensional Schr\"odinger operators.  Invent. Math.  135  (1999),
329--367.

\bibitem[K]{K} Kotani, S. Generalized Floquet theory for stationary
Schrödinger operators in one dimension.  Chaos Solitons Fractals  8  (1997),
1817--1854.

\bibitem[R]{R} Remling, C. The absolutely continuous spectrum of Jacobi
matrices.  Preprint (www.arXiv.org).

\bibitem[S1]{S2} Simon, Barry Kotani theory for one-dimensional
stochastic Jacobi matrices.  Comm. Math. Phys.  89  (1983),
227--234.

\bibitem[S2]{S3} Simon, Barry Bounded eigenfunctions and absolutely
continuous spectra for one-dimensional Schr\"odinger operators.  Proc. Amer.
Math. Soc.  124  (1996), 3361--3369.

\bibitem[S3]{S} Simon, Barry Schr\"odinger operators in the twenty-first
century.  Mathematical physics 2000,  283--288, Imp. Coll. Press, London,
2000.

\bibitem[T]{T} Trent, Tavan T. A new estimate for the vector valued corona
problem.  J. Funct. Anal.  189  (2002), 267--282.

\bibitem[U]{U} Uchiyama, A., Corona theorems for countably
many functions and estimates for their solutions, preprint, 1980.

%\bibitem[vM]{vM} van Mouche, Pierre The coexistence problem for the discrete
%Mathieu operator.  Comm. Math. Phys.  122  (1989),  no. 1, 23--33.

\end{thebibliography}
\end{document}